\documentclass[11pt,reqno]{NumPDEsArticle}

\usepackage[T1]{fontenc}
\usepackage[english]{babel}
\usepackage{csquotes}
\usepackage{amsmath,amssymb}
\usepackage{mathrsfs}
\usepackage[basicdelimiters]{NumPDEsMacros}
\usepackage{enumitem}
\usepackage{graphicx}
\usepackage[caption=false]{subfig}
\usepackage{booktabs}
\usepackage{orcidlink}
\usepackage{tikz}
\usepackage{pgfplots}

\definecolor{WildStrawberry}{rgb}{1.0, 0.26, 0.64}

\title{Robust Multigrid on Adaptive Meshes for the Integral Fractional Laplacian}
\author{Markus Faustmann\orcidlink{0009-0007-0592-8499}}
\author{Paula Hilbert\orcidlink{0009-0005-0105-1066}}
\author{Ani Miraçi\orcidlink{0000-0003-4962-9662}}
\address{Sorbonne Université, Université Paris Cité, CNRS, 
Laboratoire Jacques-Louis Lions, LJLL, F-75005 Paris, France}
\email{\tt ani.miraci@sorbonne-universite.fr}
\author{Dirk Praetorius\orcidlink{0000-0002-1977-9830}}
\address{TU Wien, Institute of Analysis and Scientific Computing, 
Wiedner Hauptstr. 8--10/E101/4, 1040 Vienna, Austria}
\email{\tt markus.faustmann@asc.tuwien.ac.at}
\email{\tt paula.hilbert@asc.tuwien.ac.at
\quad {\textrm{(corresponding author)}}}
\email{\tt dirk.praetorius@asc.tuwien.ac.at}

\dedicatory{Dedicated to Olaf Steinbach on the occasion of his 60th birthday}

\thanks{This research was funded in whole or in part by the Austrian 
Science Fund (FWF) projects
\href{https://www.fwf.ac.at/en/research-radar/10.55776/PAT3699424}{10.55776/PAT3699424}
(standalone project PAT3699424 ``Optimal robust solvers for reliable and efficient AFEMs''),
\href{https://www.fwf.ac.at/en/research-radar/10.55776/PAT3446525}{10.55776/PAT3446525}
(standalone project PAT3446525 ``Adaptive Uzawa-type FEM for nonlinear PDEs''), and 
\href{https://www.fwf.ac.at/en/research-radar/10.55776/PAT1917426}{10.55776/PAT1917426}
(standalone project PAT1917426 ``Computational fractional diffusion'').
Additionally, Paula Hilbert is supported by the Vienna School of Mathematics.}

\begin{document}

\maketitle

\begin{abstract}
We analyze a geometric multigrid method on locally refined mesh hierarchies for non-local variational problems,
focusing on the integral fractional Laplacian. For newest-vertex-bisection-generated hierarchies, we prove uniform contraction of the multigrid iteration independently of the refinement level, the smallest local mesh size, and the polynomial degree in the finite element discretization. The proof is based on an abstract Hilbert-space framework, a strengthened Cauchy--Schwarz estimate obtained from a near-/far-field decomposition of the non-local bilinear form, and 
a 
multilevel stable decomposition independent of the discretization parameters. Numerical experiments underline the theoretical findings.
\end{abstract}


\section{Introduction}


Adaptive finite element methods (AFEMs) provide an effective framework for the numerical solution of partial differential equations (PDEs), in particular for the
approximation of solutions with singular behaviour induced by either geometry, data, or the differential operator itself. Rate-optimal convergence of AFEMs is meanwhile considerably well established, see, e.g., for linear elliptic PDEs~\cite{doerfler1996, MNS00, bdd2004, stevenson2007, ckns2008} or (boundary) integral equations~\cite{FKMP2013, Gantumur2013, FMP2021}, which all fit into an abstract framework,~\cite{axioms}. 

From the outset, adaptive wavelet methods have been concerned with optimal computational complexity; see, e.g.,~\cite{cdd2001,cdd2003}.
More recently, optimal computational complexity of standard AFEMs of the form \texttt{SOLVE}-\texttt{ESTIMATE}-\texttt{MARK}-\texttt{REFINE} to steer possible local mesh refinement has come into focus.
This distinction between approximation optimality and computational optimality motivates the use of iterative algebraic solvers in the module \texttt{SOLVE} with contraction properties that are uniform along the entire adaptive mesh hierarchy. In particular, the number of algebraic iterations required to reduce the solver error should remain bounded
independently of the refinement level and the smallest local mesh size, and, for higher-order discretizations, independently of the polynomial degree. In combination with suitable stopping criteria that balance algebraic and discretization errors, uniformly contractive solvers can be incorporated into AFEM without destroying its optimal convergence properties, \cite{stevenson2007,ghps2021,bfmps2025}; see also~\cite{FHPS2019} in the context of adaptive boundary element methods.

In this work, we consider variational problems of the form 
$$
\text{find } u^\star \in \XX:\quad a(u^\star,v) = f(v) \qquad \text{for all } v \in \XX,
$$
where $a(\cdot,\cdot)$ is a symmetric, bounded, and elliptic bilinear form on a
Hilbert space $\XX$. Within this setting, our principal model problem is the integral fractional Laplacian
$$
(-\Delta)^s u=f\quad\text{in }\Omega, 
\qquad
u=0\quad\text{in }\mathbb{R}^d\setminus\Omega,
\qquad \text{for some fixed } s\in(0,1),
$$
where $(-\Delta)^s$ is understood in its singular-integral sense. We refer to \cite{Kwasnicki2017} for equivalent representations of the fractional Laplacian in $\mathbb{R}^d$ and to \cite{Lischke2020} for a comparison of different realizations of fractional Laplacians on bounded domains.
In particular, the key feature of the fractional Laplacian is that the associated energy space is a fractional-order Sobolev space, and the corresponding bilinear form $a(\cdot,\cdot)$ is non-local. 
Solving the resulting discrete systems thus requires multilevel solvers
tailored to fractional-order norms and non-local operators.

Multilevel methods as solvers or preconditioners 
on locally refined meshes for local PDEs have a rich history; we mention among many others, \cite{Hackbusch1985,BPX90,Oswald1994,cnx2012}. 
Related ideas have subsequently been developed for fractional and, more generally, non-local operators.
For the extension formulation of fractional diffusion, \cite{CNOTS2016} develops multilevel methods with nearly uniform convergence, and additive multigrid preconditioners on quasi-uniform meshes are studied in \cite{BKM2019}. For the integral fractional Laplacian, \cite{AG2017} develops an adaptive finite element framework and discusses iterative solvers, including geometric multigrid, on locally refined meshes.

Closely related topics include the stability and decomposition results in fractional Sobolev spaces. Stability of finite element projections under local mesh conditions was established in \cite{Steinbach2001}, motivated in particular by applications to preconditioning on adaptive finite and boundary element meshes. Stable decompositions for fractional Sobolev spaces and their use in
discretization-robust 
preconditioning have subsequently been studied in \cite{KMR2020,FMP2021_SZ,BNWX2023,KMR2025}.
There is also an extensive literature for boundary integral equations, where non-local operators and fractional-order Sobolev spaces arise naturally. Multilevel representations of Sobolev norms in adaptive boundary element methods were already employed in \cite{Steinbach2000}; see also~\cite{MS1999,Steinbach2003,MST00,Maischak2009,FMPR2015,FFPS2017,FHPS2019} for multilevel and (optimal) additive Schwarz-type preconditioners for boundary integral operators on (locally) refined meshes. 

For second-order linear elliptic PDEs on meshes generated by newest vertex bisection (NVB), \cite{IMPS2024} introduced a local geometric multigrid method whose contraction factor is robust with respect to both local mesh size $h$ and the polynomial degree $p$ of the finite element discretization. The method replaces the exact solve in each step of the adaptive algorithm by a V-cycle consisting of an exact solve on a fixed coarse mesh, inexpensive first-order local corrections on intermediate levels, and a higher-order vertex-patch correction on the finest level. Moreover, the multigrid iteration provides a built-in estimator for the algebraic error. The method relies on earlier work~\cite{MPV2021}, which designed and analyzed a $p$-robust multigrid method using a high-order one-level stable decomposition established in~\cite{SMPZ08} together with a lowest-order multilevel decomposition from~\cite{cnx2012}.

The non-locality of the considered operators, however, does not allow a direct transfer of the contraction argument for the multigrid solver from \cite{IMPS2024} to the non-local setting. Indeed, therein, technical requirements such as a strengthened Cauchy-Schwarz inequality and existence of a stable decomposition heavily used locality of the bilinear form and energy norm. 

The main contribution of the present work is to extend this local geometric multigrid framework to the considered non-local problems on NVB-generated adaptive hierarchies. We develop an abstract Hilbert-space framework and prove uniform contraction of the resulting multigrid iteration for the integral fractional Laplacian. The key analytical ingredients are a strengthened Cauchy--Schwarz estimate obtained through a decomposition into near- and far-field interactions 
for lowest-order discretizations and 
an $h$- and $p$-robust multilevel stable decomposition. In particular, the contraction estimate is uniform with respect to the adaptive mesh hierarchy.

We emphasize that the present paper focuses on the robust algebraic solver and the multilevel analysis required for its uniform contraction. A complete theory of optimal complexity, i.e., optimal convergence rates with respect to the total
computational cost, is beyond the scope of this article. 
We note, however, that besides the usual AFEM ingredients, it requires, in
particular for non-local problems, matrix compression techniques and an
efficient evaluation of the a posteriori error indicators.

The remainder of the paper is organized as follows. Section~\ref{sec:setting} introduces the model problem and the  fractional Sobolev setting, as well as the discretization by finite element spaces based on an NVB-generated mesh hierarchy. In section~\ref{sec:mg}, we formulate the geometric multigrid iteration in an abstract Hilbert-space framework and establish a quasi-Pythagorean estimate, reliability of a solver-inherent a-posteriori estimator for the solver error, and the resulting uniform contraction of the proposed multigrid solver. Section~\ref{sec:main} tailors the abstract theory to the fractional Laplacian, states the main result on uniform contraction for the multigrid solver and discusses multigrid as a solver for AFEM with optimal convergence rates. Sections~\ref{sec:SCS} and \ref{sec:stable_decomp} provide the essential ingredients for the proof of the uniform contraction, a strengthened Cauchy--Schwarz inequality and the existence of a 
robust multilevel stable decomposition. To conclude, Section~\ref{sec:num} provides numerical experiments highlighting the theoretical results.


\section{Setting}\label{sec:setting}


Let $\Omega \subset \mathbb{R}^d$, $d \geq 1$, be a bounded Lipschitz domain 
and let $0 < s < 1$. 
The fractional Sobolev spaces $H^s(\R^d)$ are defined, e.g., by the Fourier transform; see, e.g.,~\cite[Section~7]{lm1972} and~\cite[Definition~1.3.1.3]{Pier2011}.
We denote the distance to the boundary by 
$\operatorname{dist}(\cdot, \partial\Omega) : \Omega \to \mathbb{R}_{>0}$, 
$x \mapsto \inf_{y \in \partial\Omega} \abs{x-y}$. 
The Sobolev space $H^s(\mathbb{R}^d)$ can be equipped with the norm
\begin{equation*}
  \|u\|_{H^s(\R^d)}^2
  \coloneqq 
  \|u\|_{L^2(\mathbb{R}^d)}^2 + \abs{u}^2_{H^s(\R^d)}
  \coloneqq
   \|u\|_{L^2(\mathbb{R}^d)}^2 + \int_{\R^d} \int_{\R^d} 
  \frac{\abs{u(x) - u(y)}^2}{\abs{x-y}^{d+2s}} \d{x} \d{y}.
\end{equation*}
Moreover, we introduce the space and norm
\begin{equation}\label{eq:tHs}
	\widetilde{H}^s(\Omega) \coloneqq 
  \set{ u \in H^s(\mathbb{R}^d) \given u \equiv 0 \text{ on } \Omega^c },
	\qquad
	\|v\|_{\widetilde{H}^s(\Omega)}^2 \coloneqq \|v\|_{H^s(\mathbb{R}^d)}^2 
  + \|\operatorname{dist}(\cdot,\partial\Omega)^{-s} v\|_{L^2(\Omega)}^2.
\end{equation}
The fractional Laplacian $(-\Delta)^s$ 
can be defined via the principal value integral representation 
\begin{equation}\label{eq:fracLap}
	(-\Delta)^s u(x) 
  \coloneqq C(d,s) \operatorname{P.V.} \int_{\mathbb{R}^d} 
  \frac{u(x) - u(y)}{\abs{x-y}^{d+2s}} \, \mathrm{d}y
	\quad \text{with} \quad
	C(d,s) \coloneqq -2^{2s} \, \frac{\Gamma(s + d/2)}{\pi^{d/2} \Gamma(-s)},
\end{equation}
where $\Gamma(\,\cdot\,)$ denotes the Gamma function.
We consider the homogeneous Dirichlet problem
\begin{equation}\label{eq:strong_form}
	(-\Delta)^s u = f \quad \text{in } \Omega, 
  \qquad u = 0 \quad \text{in } \Omega^c,
\end{equation}
for a given right-hand side 
$f \in \widetilde{H}^{-s}(\Omega) \coloneqq (\widetilde{H}^s(\Omega))'$.
The weak formulation of~\eqref{eq:strong_form} reads: 
find $u^\star \in \widetilde{H}^s(\Omega)$ such that
\begin{equation}\label{eq:weak_form}
	a(u^\star, v) 
  \coloneqq \frac{C(d,s)}{2}
	\iint_{\mathbb{R}^d \times \mathbb{R}^d}
	\frac{(u^\star(x) - u^\star(y))(v(x) - v(y))}{\abs{x-y}^{d+2s}} \, 
  \mathrm{d}x \, \mathrm{d}y
   = \langle f, v \rangle 
  \;\, \text{for all } v \in \widetilde{H}^s(\Omega),
\end{equation}
Here, $\dual{\cdot}{\cdot}$ denotes the duality pairing between $\widetilde{H}^{-s}(\Omega)$ and $\widetilde{H}^s(\Omega)$, which conincides with the $L^2(\Omega)$ scalar product if $f \in L^2(\Omega)$.
The form $a(\cdot, \cdot)$ is continuous and 
coercive on the Hilbert space $\widetilde{H}^s(\Omega)$
(see, e.g., \cite{AB2017}). Therefore, existence and uniqueness of a weak solution 
$u \in \widetilde{H}^s(\Omega)$
to~\eqref{eq:weak_form} follow from the Lax--Milgram lemma.
We denote the induced energy norm by $\enorm{v}^2 \coloneqq a(v,v)$.

\subsection{Multilevel notation}\label{sec:multilevel}
Let \(\TT_0\) be a conforming triangulation of \(\Omega\) into compact simplices that is admissible in the sense of~\cite{stevenson2008}.
Consider a sequence of successively refined triangulations $(\TT_{\ell})_{\ell \in \N_0}$ by newest-vertex bisection (NVB) starting from the initial mesh $\TT_0$; see, e.g.,~\cite{traxler1997, stevenson2008} for $d \geq 1$ and~\cite{AFFKP13} for $d=1$.
There holds $\TT_{\ell} = \texttt{refine}(\TT_{\ell-1}, \MM_{\ell-1})$ for $\ell \geq 1$ and a set of marked elements $\MM_{\ell-1} \subseteq \TT_{\ell-1}$, such that $\TT_{\ell}$ is the coarsest conforming refinement of $\TT_{\ell-1}$ with $\MM_{\ell-1} \subseteq \TT_{\ell-1} \backslash \TT_{\ell}$.
Denote by $\T \coloneqq \T(\TT_0)$ the set of all meshes that can be generated via NVB from $\TT_0$. 
Meshes generated by NVB are uniformly $\gamma$-shape regular, i.e., 
\begin{equation*}
  \max_{\TT \in \T}\max_{T \in \TT} \frac{\diam(T)}{\abs{T}^{1/d}} \leq \gamma < \infty,
\end{equation*}
where $\gamma > 0$ depends only on $\TT_0$; see, e.g.,~\cite[Theorem 2.1]{stevenson2008} for $d \geq 2$ and~\cite{AFFKP13} for $d=1$.
For all subsets $\omega \subseteq \Omega$, we define the $n$-patch $\TT_{\ell}^n(\omega)$ inductively by
\begin{equation*}
	\TT_{\ell}(\omega) \coloneqq \TT_{\ell}^1(\omega) 
  \coloneqq \set{T \in \TT_{\ell} \given \omega \cap T \neq \emptyset}, \qquad
	\TT_{\ell}^{n+1}(\omega) 
  \coloneqq \set{T \in \TT_{\ell} \given \exists T' \in \TT_{\ell}^n(\omega) \text{ with } T \cap T' \neq \emptyset}.
\end{equation*}
The corresponding $n$-patch subdomains are denoted by
\begin{equation*}
	\omega_{\ell}^n(\omega) \coloneqq \operatorname{interior}\Bigl(\bigcup_{T \in \TT_{\ell}^n(\omega)} T\Bigr).
\end{equation*}
Let $\VV_{\ell}$ denote the set of vertices of $\TT_{\ell}$.
For $z \in \VV_{\ell}$, we write $\TT^n_{\ell}(z) \coloneqq \TT^n_{\ell}(\{z\})$ as well as $\omega_{\ell}^n(z) \coloneqq \omega_\ell^n(\{z\})$.
Lastly, the size of a vertex-patch subdomain is given by
\begin{equation*}
	h_{\ell,z} \coloneqq \max_{T \in \TT_{\ell}(z)} h_T  \qquad \text{with} \quad h_T \coloneqq \abs{T}^{1/d}.
\end{equation*}
We define the set of new vertices of $\VV_{\ell}$ and their direct neighbors in $\VV_{\ell} \cap \VV_{\ell-1}$ by
\begin{equation}\label{def:new_vertices}
	\VV_0^+ \coloneqq \VV_0 \quad \text{and} \quad
	\VV_{\ell}^+ \coloneqq \VV_{\ell} \backslash \VV_{\ell-1} \cup \set{z \in \VV_{\ell} \cap \VV_{\ell-1} \given \omega_{\ell}(z) \neq \omega_{\ell-1}(z)} \quad \text{for } \ell \geq 1.
\end{equation}
Given a triangulation $\TT_{\ell}$ and a polynomial degree $p \in \N$, we define the finite element space
\begin{equation*}
  \XX_{\ell}^p \coloneqq \set{v \in H^1_0(\Omega) \given v|_T \in \mathbb{P}^p(T) \text{ for all } T \in \TT_{\ell}},
\end{equation*}
where $\mathbb{P}^p(T)$ denotes the space of polynomials of degree at most $p$ on the simplex $T$.
The discrete problem associated with~\eqref{eq:weak_form} reads: find $u_\ell^\star \in \XX_\ell^p$ such that
\begin{equation}\label{eq:discrete_problem}
  a(u_\ell^\star, v_\ell) = \langle f, v_\ell \rangle \qquad \text{for all } v_\ell \in \XX_\ell^p.
\end{equation}
There holds nestedness of the discrete spaces in the sense that
\begin{equation*}
  \XX_0^1 \subseteq \XX_1^1 \subseteq \XX_2^1 \subseteq \cdots \subseteq \XX_{\ell-1}^1 \subseteq \XX_{\ell}^p \subseteq \widetilde{H}^s(\Omega).
\end{equation*}
Moreover, we define the local subspaces for all $z \in \VV_{\ell}$ and $q \in \{1,p\}$ by
\begin{equation*}
  \XX_{\ell,z}^q \coloneqq \set{v_{\ell} \in \XX_{\ell}^q \given v_{\ell}|_T = 0 \text{ for all } T \in \TT_{\ell} \backslash \TT_{\ell}(z)}.
\end{equation*}

Consider the sequence of uniformly refined triangulations $(\widehat{\mathcal{T}}_{\ell})_{\ell \in \mathbb{N}_0}$, defined by $\widehat{\mathcal{T}}_0 \coloneqq \mathcal{T}_0$ and $\widehat{\mathcal{T}}_{\ell} \coloneqq \texttt{refine}(\widehat{\mathcal{T}}_{\ell-1}, \widehat{\mathcal{T}}_{\ell-1})$ for $\ell \geq 1$, whose mesh-size satisfies $\widehat{h}_{\ell} \coloneqq \max_{T \in \widehat{\mathcal{T}}_{\ell}} |T|^{1/d} = 2^{-\ell/d} \widehat{h}_0$.
Since the initial mesh is admissible and NVB conserves admissibility, every element $T \in \widehat{\TT}_{\ell-1}$ is bisected only once into two children $T_1, T_2 \in \widehat{\TT}_{\ell}$ with $\abs{T_1} = \abs{T_2} = \abs{T}/2$; see~\cite{stevenson2008} for more details. 
We will denote any object associated with the quasi-uniform mesh $\widehat{\TT}_{\ell}$ by a hat, e.g., $\widehat{\VV}_{\ell}$ and $\widehat{\omega}_{\ell}^n(z)$.


\section{Optimal multigrid solver}\label{sec:mg}


In this section, we present the geometric multigrid solver in an abstract framework and prove uniform contraction under suitable assumptions. 
Let $\XX$ be a Hilbert space equipped with a scalar product $a(\cdot,\cdot)$ that is continuous and coercive on $\XX$ and induces the energy norm $\enorm{v}^2 \coloneqq a(v,v)$.
Let $F\in \XX'$ be arbitrary.
We consider the following variational problem: find $u^\star \in \XX$ such that
\begin{equation*}
  a(u^\star,v) = F(v) \qquad \text{for all } v \in \XX.
\end{equation*}
Let $\ell \in \N_0$ and $\XX_{\ell} \subset \XX$ be a finite-dimensional subspace of $\XX$ and let $u_{\ell}^\star \in \XX_{\ell}$ be the Galerkin approximation of $u^\star$, i.e.,
\begin{equation}\label{eq:galerkin_problem}
  a(u_{\ell}^\star,v_{\ell}) = F(v_{\ell}) \qquad \text{for all } v_{\ell} \in \XX_{\ell}.
\end{equation}
The continuous problem as well as the discrete one is well-posed by the Lax-Milgram lemma.
Suppose that $\XX_{\ell}$ admits a multilevel decomposition into subspaces $\XX_{\ell',z} \subset \XX_{\ell}$, indexed by levels $\ell' \in \{0, \dots, \ell\}$ and, for each level $\ell'$, by $z \in \ZZ_{\ell'}$, where $\ZZ_{\ell'}$ is a finite index set, such that
\begin{equation*}
  \XX_{\ell} = \sum_{\ell'=0}^\ell \sum_{z \in \ZZ_{\ell'}} \XX_{\ell',z}.
\end{equation*}

\subsection{Multigrid algorithm with local subspace corrections}

Consider the following multigrid algorithm for the solution of the Galerkin problem~\eqref{eq:galerkin_problem} on level $\ell$.
Given a current approximation $u_\ell \in \XX_\ell$ of $u_\ell^\star$, define the residual functional $R_\ell \colon \XX_\ell \to \R$ by $R_\ell(v_\ell) \coloneqq F(v_\ell) - a(u_\ell, v_\ell) = a(u_\ell^\star - u_\ell, v_\ell)$ for all $v_\ell \in \XX_\ell$.

\begin{algorithm}[V-cycle algorithm]\label{alg:multigrid}
\textbf{Input:} Current approximation $u_\ell \in \XX_\ell$, space decomposition $\XX_{\ell',z}$ for $\ell' \in \{0, \dots, \ell\}$ and $z \in \ZZ_{\ell'}$. Initialize $\zeta_{-1} \coloneqq 0$ and $\sigma_{-1} \coloneqq 0$. \\
\noindent For all $\ell'=0,\dots,\ell$, perform the following steps \textnormal{(i)--(ii)}:
\begin{enumerate}[label=\textnormal{(\roman*)}, ref = \textnormal{\roman*}]
  \item \label{algo:mg:step1} For every $z \in \ZZ_{\ell'}$, define the local residual lifting $\rho_{\ell',z} \in \XX_{\ell',z}$ by
  \begin{equation}\label{eq:def_local_updates}
    a(\rho_{\ell',z}, v_{\ell',z}) = R_\ell(v_{\ell',z}) - a(\sigma_{\ell'-1}, v_{\ell',z}) \qquad \text{for all } v_{\ell',z} \in \XX_{\ell',z}.
  \end{equation}
  \item \label{algo:mg:step2} Choose a step-size $\lambda_{\ell'} > 0$ and define the level update $\rho_{\ell'} \coloneqq \sum_{z \in \ZZ_{\ell'}} \rho_{\ell',z} \in \XX_{\ell'}$ as well as 
  \begin{equation}\label{eq:def_global_updates}
    \sigma_{\ell'} \coloneqq \sigma_{\ell'-1} + \lambda_{\ell'} \rho_{\ell'} \qquad \text{and} \qquad \zeta_{\ell'}^2 \coloneqq \zeta_{\ell'-1}^2 + \lambda_{\ell'} \sum_{z \in \ZZ_{\ell'}} \enorm{\rho_{\ell',z}}^2.
  \end{equation}
\end{enumerate}
\textbf{Output:} Updated approximation $u_\ell + \sigma_\ell \in \XX_\ell$ and a-posteriori estimator $\zeta_\ell(u_\ell) \coloneqq \zeta_\ell$.
\end{algorithm}

\subsection{Uniform contraction}

Under the following assumptions on the space decomposition and the step-sizes, one iteration of the V-cycle (see Algorithm~\ref{alg:multigrid}) contracts the algebraic error uniformly in the energy norm; see Theorem~\ref{theorem:uniform_contraction} below:

\begin{assumption}[choice of step-size]\label{assumption:one}
For the choice of the step-sizes $\lambda_{\ell'}$, we suppose the following:
\begin{enumerate}[label=\textnormal{(\alph*)}, ref = \textnormal{\alph*}]
  \item \label{assumption:one:a} uniform bounds: There exist constants $0 < \lambda_\ast < \lambda^\ast$ independent of $\ell$ such that
\begin{equation*}
\lambda_\ast \le \lambda_{\ell'}
\quad\text{for all } \ell'= 0, \dots, \ell \qquad \text{and} \qquad \lambda_{\ell'} \le \lambda^\ast \quad\text{for all } \ell'= 0, \dots, \ell-1.
\end{equation*}
\item \label{assumption:one:b} levelwise update bound: For every level $\ell' \in \{0,\dots,\ell\}$, there holds
\begin{equation}\label{eq:assumption2}
\enorm{\lambda_{\ell'}\rho_{\ell'}}^2
\le
\lambda_{\ell'} \sum_{z\in\ZZ_{\ell'}}\enorm{\rho_{\ell',z}}^2.
\end{equation}
\end{enumerate}
\end{assumption}

\begin{assumption}[strengthened Cauchy--Schwarz inequality]\label{assumption:two}
There exists a constant $\Cscs > 0$, independent of $\dim \XX_\ell$, such that for all families $v_{\ell',z}, u_{\ell',z} \in \XX_{\ell',z}$ with $\ell' \in \{0, \dots, \ell\}$ and $z \in \ZZ_{\ell'}$, there holds
\begin{equation}\label{eq:assumption3}
\sum_{\ell'=0}^{\ell} \sum_{k=0}^{\ell'-1}\sum_{z\in\ZZ_{\ell'}}\sum_{w\in\ZZ_{k}} a(v_{\ell',z}, u_{k,w})
\le
\Cscs
\Big( \sum_{\ell'=0}^\ell\sum_{z\in\ZZ_{\ell'}}\enorm{v_{\ell',z}}^2 \Big)^{1/2} \Big(\sum_{k=0}^{\ell-1}\sum_{w\in\ZZ_k}\enorm{u_{k,w}}^2 \Big)^{1/2}.
\end{equation}
\end{assumption}

\begin{assumption}[stable decomposition]\label{assumption:three}
There exists a constant $\Csd \ge 1$, independent of $\dim \XX_\ell$, such
that every $v_{\ell} \in \XX_\ell$ admits a stable decomposition satisfying
\begin{equation*}
  v_\ell
  =
  \sum_{\ell'=0}^{\ell} \sum_{z \in \ZZ_{\ell'}} v_{\ell',z},
  \quad
  v_{\ell',z} \in \XX_{\ell',z},
\qquad
\text{and}
\qquad
\sum_{\ell'=0}^{\ell}\sum_{z\in\ZZ_{\ell'}}\enorm{v_{\ell',z}}^2
\le
\Csd \enorm{v_\ell}^2.
\end{equation*}
\end{assumption}

\begin{remark}
Assumptions~\ref{assumption:two}--\ref{assumption:three} are standard in the literature on multilevel methods; see, e.g.,\linebreak \cite[Assumption~2.2~\&~2.3]{TW05} and~\cite[Assumption~(A1)~\&~(A2)]{cnx2012}. Usually, the literature also imposes an additional assumption on the levelwise operators; see, e.g.,~\cite[Section~2]{BPWX91} or~\cite[Assumption~(T)]{cnx2012}. 
This corresponds to Assumption~\ref{assumption:one}\eqref{assumption:one:b}.
Indeed, define the levelwise operator $\SS_{\ell'} \colon \XX_\ell \to \XX_{\ell'}$
by $\SS_{\ell'}(v_\ell) \coloneqq \lambda_{\ell'} \sum_{z \in \ZZ_{\ell'}}
v_{\ell',z}$, where 
\begin{equation*}
  a(v_{\ell',z}, w_{\ell',z}) = a(v_\ell, w_{\ell',z})
  \ \text{for all } w_{\ell',z} \in \XX_{\ell',z}
  \quad \text{and} \quad 
  \lambda_{\ell'} \, \enorm{\sum_{z \in \ZZ_{\ell'}} v_{\ell',z}}^2 \le \sum_{z \in \ZZ_{\ell'}} \enorm{v_{\ell',z}}^2.
\end{equation*}
Here, the estimate on the step-size is the analog of
Assumption~\ref{assumption:one}\eqref{assumption:one:b} for arbitrary $v_\ell \in \XX_\ell$. 
Note that Assumption~\ref{assumption:one}\eqref{assumption:one:b} is less restrictive, in the sense that it
requires it only for the current error
$v_\ell = u_\ell^\star - u_\ell - \sigma_{\ell'-1}$, where
$v_{\ell',z} = \rho_{\ell',z}$. Nonetheless, we can show that $\SS_{\ell'}$ is stable with constant $1$, i.e., $\enorm{\SS_{\ell'}(v_\ell)} \le \enorm{v_\ell}$ for all $v_\ell \in \XX_\ell$.
To see this (in the case that both $v_\ell\neq 0$ and $\sum_{z \in \ZZ_{\ell'}}
v_{\ell',z}\neq 0$, while the claim holds trivially otherwise), note that using the definition of  $\SS_{\ell'}$, the estimate on the step-size, and the local projection property, we obtain
\begin{equation*}
  \enorm{\SS_{\ell'}(v_\ell)}
  =
  \enorm{\lambda_{\ell'} \sum_{z \in \ZZ_{\ell'}} v_{\ell',z}}
   \le \frac{\sum_{z \in \ZZ_{\ell'}} \enorm{v_{\ell',z}}^2}{\enorm{\sum_{z \in \ZZ_{\ell'}} v_{\ell',z}}}
   = \frac{\sum_{z \in \ZZ_{\ell'}} 
   a(v_{\ell},v_{\ell',z})}{\enorm{\sum_{z \in \ZZ_{\ell'}} v_{\ell',z}}}
   = \frac{ 
   a(v_{\ell},\sum_{z \in \ZZ_{\ell'}} v_{\ell',z})}{\enorm{\sum_{z \in \ZZ_{\ell'}} v_{\ell',z}}}
   \le   \enorm{v_\ell}.
\end{equation*}
%
%
\end{remark}

Let $\Phi_\ell \colon \XX_\ell \to \XX_{\ell}$ denote the iteration operator of the V-cycle, i.e., $\Phi_\ell(u_\ell) \coloneqq u_\ell + \sigma_\ell$ and $\sigma_\ell$ is defined in Algorithm~\ref{alg:multigrid}.
Having 
identified the above assumptions, we can now state the main result of this section, which is an abstraction of~\cite[ Theorem~2.5]{IMPS2024}.

\begin{theorem}[uniform contraction of the multigrid]
\label{theorem:uniform_contraction}
Suppose that Assumptions~\ref{assumption:one}--\ref{assumption:three} hold and let $u_\ell^\star \in \XX_\ell$ denote the solution of~\eqref{eq:galerkin_problem}. 
Then, the V-cycle iteration from Algorithm~\ref{alg:multigrid} satisfies the following quasi-Pythagorean estimate
\begin{equation}\label{eq:pythagoras}
  \enorm{u_\ell^\star - \Phi_\ell(u_\ell)}^2
  \leq 
  \enorm{u_\ell^\star - u_\ell}^2 - \zeta_\ell(u_\ell)^2
  \qquad \text{for all } u_\ell \in \XX_\ell.
\end{equation}
Moreover, the built-in a-posteriori estimator $\zeta_\ell(u_\ell)$ is reliable, i.e., there holds 
\begin{equation}\label{eq:reliability}
  \enorm{u_\ell^\star - u_\ell} \leq \Crel \, \zeta_\ell(u_\ell) \qquad \text{for all } u_\ell \in \XX_\ell,
\end{equation}
where $\Crel \coloneqq \sqrt{\Csd(1+ \lambda_\ast^{-1})(1+ \Cscs^2\lambda^\ast)}$.
In particular, these two estimates imply uniform contraction
\begin{equation}\label{eq:contraction}
  \enorm{u_\ell^\star - \Phi_\ell(u_\ell)} \le \qmg \, \enorm{u_\ell^\star - u_\ell} \quad \text{for all } u_\ell \in \XX_\ell \quad \text{with} \quad 0 < \qmg \coloneqq \sqrt{1-\Crel^{-2}} < 1.
\end{equation}
\end{theorem}

\begin{proof}
The proof is split into three steps.

\textbf{Step 1 (quasi-Pythagorean estimate):}
Let $u_\ell \in \XX_\ell$ and define 
\begin{equation*}
  e_{-1} \coloneqq u_\ell^\star - u_\ell, \quad  e_{\ell'} \coloneqq u_\ell^\star - u_\ell - \sigma_{\ell'} \quad \text{for all } \ell' \in \{0, \dots, \ell\}.
\end{equation*}
Since $e_{\ell'} = e_{\ell'-1} - \lambda_{\ell'} \rho_{\ell'}$, Assumption~\ref{assumption:one}\eqref{assumption:one:b} implies that
\begin{equation*}
  \enorm{e_{\ell'}}^2
  =
  \enorm{e_{\ell'-1}}^2 - 2 a(e_{\ell'-1}, \lambda_{\ell'} \rho_{\ell'}) + \enorm{\lambda_{\ell'} \rho_{\ell'}}^2
  \eqreff{eq:assumption2}{\leq}
  \enorm{e_{\ell'-1}}^2 - 2 a(e_{\ell'-1}, \lambda_{\ell'} \rho_{\ell'}) + \lambda_{\ell'} \sum_{z \in \ZZ_{\ell'}} \enorm{\rho_{\ell',z}}^2.
\end{equation*}
Due to the definition~\eqref{eq:def_local_updates} of $\rho_{\ell',z}$, there holds
\begin{equation*}
  a(e_{\ell'-1}, \lambda_{\ell'} \rho_{\ell'})
  =
  \lambda_{\ell'} \sum_{z \in \ZZ_{\ell'}} a(e_{\ell'-1}, \rho_{\ell',z})
  =
  \lambda_{\ell'} \sum_{z \in \ZZ_{\ell'}} \big(R_\ell(\rho_{\ell',z}) - a(\sigma_{\ell'-1}, \rho_{\ell',z})\big)
  \eqreff{eq:def_local_updates}{=}
  \lambda_{\ell'} \sum_{z \in \ZZ_{\ell'}} \enorm{\rho_{\ell',z}}^2.
\end{equation*}
Hence, we obtain
\begin{equation*}
  \enorm{e_{\ell'}}^2
  \leq
  \enorm{e_{\ell'-1}}^2 - \lambda_{\ell'}
\sum_{z \in \ZZ_{\ell'}} \enorm{\rho_{\ell',z}}^2 \quad \text{for all } \ell' \in \{0, \dots, \ell\}.
\end{equation*}
By induction, we conclude that
\begin{equation*}
  \enorm{u_\ell^\star - \Phi_\ell(u_\ell)}^2
  =
  \enorm{e_\ell}^2
  \leq
  \enorm{e_{-1}}^2 - \sum_{\ell'=0}^{\ell} \lambda_{\ell'} \sum_{z \in \ZZ_{\ell'}} \enorm{\rho_{\ell',z}}^2
  \eqreff{eq:def_global_updates}{=}
  \enorm{u_\ell^\star - u_\ell}^2 - \zeta_\ell(u_\ell)^2.
\end{equation*}

\textbf{Step 2 (reliability):}
By Assumption~\ref{assumption:three}, there exists a decomposition 
\begin{equation}\label{eq:error_decomposition}
  u_\ell^\star - u_\ell = \sum_{\ell'=0}^{\ell} \sum_{z \in \ZZ_{\ell'}} u_{\ell',z} 
  \quad \text{with} \quad  
  u_{\ell',z} \in \XX_{\ell',z} 
  \quad \text{and} \quad 
  \sum_{\ell'=0}^{\ell} \sum_{z \in \ZZ_{\ell'}} \enorm{u_{\ell',z}}^2
  \leq
  \Csd \enorm{u_\ell^\star - u_\ell}^2.
\end{equation}
Using this decomposition, the definition~\eqref{eq:def_local_updates} of $\rho_{\ell',z}$, and the definition~\eqref{eq:def_global_updates} of $\sigma_{\ell'}$, we obtain
\begin{align}\label{eq:error_splitting}
\begin{split}
  &\enorm{u_\ell^\star - u_\ell}^2
  =
  a(u_\ell^\star - u_\ell, u_\ell^\star - u_\ell)
  \eqreff{eq:error_decomposition}{=}
  \sum_{\ell'=0}^{\ell} \sum_{z \in \ZZ_{\ell'}} a(u_\ell^\star - u_\ell, u_{\ell',z}) 
  =
  \sum_{\ell'=0}^{\ell} \sum_{z \in \ZZ_{\ell'}}  R_{\ell}(u_{\ell',z})\\
 &\eqreff{eq:def_local_updates}{=}
  \sum_{\ell'=0}^{\ell} \sum_{z \in \ZZ_{\ell'}} \!\! \big(a(\rho_{\ell',z}, u_{\ell',z}) \!+ \! a(\sigma_{\ell'-1}, u_{\ell',z}) \big)
  \eqreff{eq:def_global_updates}{=} 
  \sum_{\ell'=0}^{\ell} \sum_{z \in \ZZ_{\ell'}} \!\! \Big(a(\rho_{\ell',z}, u_{\ell',z}) + \sum_{k=0}^{\ell'-1} a(\lambda_k \rho_k, u_{\ell',z})\Big).
\end{split}
\end{align}
The Cauchy--Schwarz inequality and Young's inequality for arbitrary $\varepsilon > 0$ yield that
\begin{equation*}
  a(\rho_{\ell',z}, u_{\ell',z}) 
  \leq
  \enorm{\rho_{\ell',z}} \enorm{u_{\ell',z}}
  =
  (\lambda_{\ell'}^{1/2}\enorm{\rho_{\ell',z}})(\lambda_{\ell'}^{-1/2}\enorm{u_{\ell',z}}) 
  \leq 
  \frac{1}{4\varepsilon} \lambda_{\ell'}\enorm{\rho_{\ell',z}}^2 + \varepsilon \lambda_{\ell'}^{-1} \enorm{u_{\ell',z}}^2.
\end{equation*}
With $\lambda_{\ell'}^{-1} \leq \lambda_\ast^{-1}$ from Assumption~\ref{assumption:one}\eqref{assumption:one:a} and the definition~\eqref{eq:def_global_updates} of $\zeta_{\ell}$, it follows that
\begin{equation}\label{eq:nolevel_interactions}
  \sum_{\ell'=0}^{\ell} \sum_{z \in \ZZ_{\ell'}} a(\rho_{\ell',z}, u_{\ell',z}) 
  \leq
  \frac{1}{4\varepsilon} \zeta_\ell(u_\ell)^2 + \varepsilon \lambda_\ast^{-1} \sum_{\ell'=0}^{\ell} \sum_{z \in \ZZ_{\ell'}} \enorm{u_{\ell',z}}^2.
\end{equation}
For the second term in \eqref{eq:error_splitting}, we recall $\rho_k = \sum_{z \in \ZZ_k} \rho_{k,z}$ and apply Assumption~\ref{assumption:two} to obtain
\begin{equation*}
  \sum_{\ell'=0}^{\ell}\sum_{k=0}^{\ell'-1} \sum_{z \in \ZZ_{\ell'}} a(\lambda_k \rho_k, u_{\ell',z}) 
  \eqreff{eq:assumption3}{\leq}
  \Cscs
  \Big(\sum_{k=0}^{\ell-1}\sum_{z\in\ZZ_k}\lambda_k^2 \enorm{\rho_{k,z}}^2 \Big)^{1/2}
  \Big( \sum_{\ell'=0}^\ell\sum_{z\in\ZZ_{\ell'}}\enorm{u_{\ell',z}}^2 \Big)^{1/2}.
\end{equation*}
Assumption~\ref{assumption:one}\eqref{assumption:one:a} and the definition~\eqref{eq:def_global_updates} of $\zeta_\ell$ give
\begin{equation*}
  \sum_{k=0}^{\ell-1}\sum_{z\in\ZZ_k}\lambda_k^2 \enorm{\rho_{k,z}}^2
  \leq
  \lambda^\ast \sum_{k=0}^{\ell}\sum_{z\in\ZZ_k} \lambda_k \enorm{\rho_{k,z}}^2
  =
  \lambda^\ast \zeta_\ell(u_\ell)^2.
\end{equation*}
Another application of Young's inequality for arbitrary $\varepsilon > 0$ yields
\begin{equation}\label{eq:level_interactions}
  \sum_{\ell'=0}^{\ell}\sum_{k=0}^{\ell'-1} \sum_{z \in \ZZ_{\ell'}} a(\lambda_k \rho_k, u_{\ell',z}) 
  \leq 
  \frac{\Cscs^2 \lambda^\ast}{4 \varepsilon} \zeta_\ell(u_\ell)^2 + \varepsilon \sum_{\ell'=0}^{\ell} \sum_{z \in \ZZ_{\ell'}} \enorm{u_{\ell',z}}^2.
\end{equation}
Plugging~\eqref{eq:nolevel_interactions} and~\eqref{eq:level_interactions} into~\eqref{eq:error_splitting} and using stability of the decomposition~\eqref{eq:error_decomposition} proves
\begin{equation*}
  \enorm{u_\ell^\star - u_\ell}^2
  \leq
  \Big(\frac{1+\Cscs^2 \lambda^\ast}{4 \varepsilon}\Big) \,\zeta_\ell(u_\ell)^2 + \varepsilon \, (1 + \lambda_\ast^{-1}) \, \Csd \, \enorm{u_\ell^\star - u_\ell}^2. 
\end{equation*}
Choosing $\varepsilon \coloneqq 1/(2(1 + \lambda_\ast^{-1}) \Csd)$ gives
\begin{equation*}
  \enorm{u_\ell^\star - u_\ell}^2
  \leq  
  \Csd (1 + \lambda_\ast^{-1}) (1 + \Cscs^2 \lambda^\ast) \, \zeta_\ell(u_\ell)^2
\end{equation*}
and concludes the proof of reliability.

\textbf{Step 3 (uniform contraction):}
Uniform contraction follows directly from the quasi-Pythagorean estimate and the reliability estimate, i.e.,
\begin{equation*} 
  \enorm{u_\ell^\star - \Phi_\ell(u_\ell)}^2
  \eqreff{eq:pythagoras}{\leq}
  \enorm{u_\ell^\star - u_\ell}^2 - \zeta_\ell(u_\ell)^2
  \eqreff{eq:reliability}{\leq}
  (1-\Crel^{-2}) \, \enorm{u_\ell^\star - u_\ell}^2.
\end{equation*}
This concludes the proof of~\eqref{eq:contraction} with $\qmg \coloneqq \sqrt{1-\Crel^{-2}} < 1$.
\end{proof} 


\section{Main result}\label{sec:main}


In this section, we transfer the abstract theory from section~\ref{sec:mg} to the fractional Laplacian from section~\ref{sec:setting} and derive uniform contraction of the multigrid solver for the discrete problem~\eqref{eq:discrete_problem}.

Using the notation from section~\ref{sec:setting}, we choose $\XX \coloneqq \widetilde{H}^s(\Omega)$, $\XX_\ell \coloneqq \XX_\ell^p$, $\ZZ_0 \coloneqq \{0\}$, $\ZZ_{\ell'} \coloneqq \VV_{\ell'}^+$ for $\ell' \in \{1, \dots, \ell-1\}$, and $\ZZ_{\ell} \coloneqq \VV_\ell$ if $p>1$ or $\ZZ_{\ell} \coloneqq \VV_\ell^+$ if $p=1$. 
Moreover, we define $\XX_{0,0} \coloneqq \XX_0^1$, the local subspaces $\XX_{\ell',z} \coloneqq \XX_{\ell',z}^1$ for all $\ell' \in \{1, \dots, \ell-1\}$ and $z \in \VV_{\ell'}^+$, and $\XX_{\ell,z}\coloneqq \XX_{\ell,z}^p$ for all $z \in \VV_\ell$ if $p>1$ or $\XX_{\ell,z}\coloneqq \XX_{\ell,z}^1$ for all $z \in \VV_\ell^+$ if $p=1$.
We take the scalar product $a(\cdot,\cdot)$ and the linear form $F$ from~\eqref{eq:weak_form}.
Define the optimal descent levelwise step-size via line-search as in~\cite{MPV2021} for $\rho_{\ell'} \neq 0$ by
\begin{equation} \label{eq:optimal_step_size}
  s_{\ell'} \coloneqq \frac{R_{\ell}(\rho_{\ell'})-a(\sigma_{\ell'-1},\rho_{\ell'})}{\enorm{\rho_{\ell'}}^2}.
\end{equation}
Lastly, we set
\begin{equation}\label{eq:step_sizes}
  \lambda_0 \coloneqq 1, 
  \quad 
  \lambda_{\ell'} \coloneqq \begin{cases}
    s_{\ell'} & \text{if } s_{\ell'} \leq d+1, \\
    (2(d+1))^{-1} & \text{otherwise}.
  \end{cases}
  \quad \text{for all } \ell' \in \{1, \dots, \ell-1\}, 
  \quad 
  \lambda_\ell \coloneqq s_\ell.
\end{equation}
This means that on the initial mesh we use the full space $\XX_0^1$ and perform an exact solve, on the intermediate levels the local spaces are one-dimensional so that no matrix inversion is required (corresponding to Jacobi smoothing), and on the finest level the local spaces scale with the polynomial degree $p$ (corresponding to block Jacobi smoothing).
Moreover, the step-sizes are needed to compensate for the patch overlap in the multilevel decomposition, where the optimal step-size is chosen whenever possible, in the sense that $s_{\ell'}$ can be uniformly bounded, and a fixed step-size is used otherwise.
With these choices, we can show that indeed Assumptions~\ref{assumption:one}--\ref{assumption:three} from section~\ref{sec:mg} are satisfied.
For easier notation, we set $\rho_{\ell,z} \coloneqq 0$ for all $z \in \VV_\ell \setminus \VV_\ell^+$ if $p=1$.

\begin{lemma}[step-size properties]\label{lemma:step_size}
There holds 
\begin{equation*}
  \lambda_\ast \coloneqq \frac{1}{2(d+1)} \leq \lambda_{\ell'}
  \quad \text{for all } \ell' \in \set{0, \dots, \ell} 
  \quad \text{and} \quad 
  \lambda_{\ell'} \leq d+1 \eqqcolon \lambda^\ast \quad \text{for all } \ell' \in \set{0, \dots, \ell-1}.
\end{equation*}
Moreover, we have 
\begin{equation*}
  \enorm{\lambda_{\ell'} \rho_{\ell'}}^2
  \leq
  \lambda_{\ell'} \sum_{z\in \VV_{\ell'}^+}\enorm{\rho_{\ell',z}}^2
  \qquad \text{for all } \ell' \in \set{1, \dots, \ell-1} 
  \quad \text{and} \quad 
  \enorm{\lambda_{\ell} \rho_{\ell}}^2
  \leq
  \lambda_{\ell} \sum_{z\in \VV_{\ell}}\enorm{\rho_{\ell,z}}^2.
\end{equation*}
\end{lemma}

\begin{proof}
The proof is split into three steps.

\textbf{Step 1 (finite patch overlap):}
Let $v_{\ell',z} \in \XX_{\ell',z}^q$ for $q \in \{1,p\}$, $z \in \VV_{\ell'}$, and some $\ell' \in \N_0$.
Fix $x,y \in \R^d$ and define $\alpha_{z}(x,y) \coloneqq v_{\ell',z}(x) - v_{\ell',z}(y)$ and $\beta_{x,y} \coloneqq \set{z \in \VV_{\ell'} \given \alpha_z(x,y) \neq 0}$.
Since $v_{\ell',z}(x) \neq 0$ only for the vertices $z$ of one element $T \in \TT_{\ell'}$ with $x \in T$, there holds $\# \beta_{x,y} \leq 2(d+1)$.
Hence, it follows from the discrete Cauchy--Schwarz inequality that
\begin{equation*}
  \Big(\sum_{z\in \VV_{\ell'}}\big(v_{\ell',z}(x)-v_{\ell',z}(y)\big)\Big)^2
  = 
  \Big(\!\! \sum_{z \in \beta_{x,y}} \!\! \alpha_z(x,y)\Big)^2
  \leq
  2(d+1) \sum_{z \in \VV_{\ell'}} \big(v_{\ell',z}(x)-v_{\ell',z}(y)\big)^2.
\end{equation*}
Dividing by $\abs{x-y}^{d+2s}$, integrating over $\R^d \times \R^d$, and multiplying by $C(d,s)/2$, yields
\begin{equation}\label{eq:patch_overlap}
  \Big\lvert \hspace{-0.2em} \Big\lvert \hspace{-0.2em} \Big\lvert \sum_{z \in \VV_{\ell'}}v_{\ell',z} \Big\rvert \hspace{-0.2em} \Big\rvert \hspace{-0.2em} \Big\rvert^2
  \leq
  2(d+1) \sum_{z \in \VV_{\ell'}} \enorm{v_{\ell',z}}^2.
\end{equation}

\textbf{Step 2 (step-size bounds):}
The upper bound follows directly from the definition~\eqref{eq:step_sizes} of $\lambda_{\ell'}$.
In the case where $s_{\ell'} > d+1$ the lower bound is also given by definition. 
Otherwise, we get by~\eqref{eq:def_local_updates} and step 1 that
\begin{equation*}
  s_{\ell'} 
  = 
  \frac{R_{\ell}(\rho_{\ell'})-a(\sigma_{\ell'-1},\rho_{\ell'})}{\enorm{\rho_{\ell'}}^2}
  \eqreff{eq:def_local_updates}{=}
  \frac{\sum_{z \in \VV_{\ell'}^+} \enorm{\rho_{\ell',z}}^2}{\enorm{\rho_{\ell'}}^2} 
  \eqreff{eq:patch_overlap}{\geq} 
  \frac{1}{2(d+1)}.
\end{equation*}
Note that we have $ z \in \VV_{\ell}$ for the finest level $\ell$ if $p > 1$, but the same arguments apply.

\textbf{Step 3 (levelwise update bound):}
Let $\lambda_{\ell'} = s_{\ell'}$. 
Due to Algorithm~\ref{alg:multigrid}\eqref{algo:mg:step1}, there holds
\begin{equation*}
  \enorm{\lambda_{\ell'}\rho_{\ell'}}^2 
  \eqreff{eq:optimal_step_size}{=} 
  \lambda_{\ell'}\frac{R_{\ell}(\rho_{\ell'})-a(\sigma_{\ell'-1},\rho_{\ell'})}{\enorm{\rho_{\ell'}}^2} \enorm{\rho_{\ell'}}^2 
  \eqreff{eq:def_local_updates}{=}
  \lambda_{\ell'} \sum_{z \in \VV_{\ell'}^+} \enorm{\rho_{\ell',z}}^2,
\end{equation*}
where again we have $z \in \VV_{\ell}$ if $p > 1$.
In the other case, finite patch overlap yields
\begin{equation*}
  \enorm{\lambda_{\ell'}\rho_{\ell'}}^2 
  \eqreff{eq:step_sizes}{=}
  \frac{\lambda_{\ell'}}{2(d+1)} \enorm{\rho_{\ell'}}^2
  \eqreff{eq:patch_overlap}{\leq}
  \frac{\lambda_{\ell'}}{2 (d+1)} 2(d+1) \sum_{z \in \VV_{\ell'}^+} \enorm{\rho_{\ell',z}}^2
  =
  \lambda_{\ell'} \sum_{z \in \VV_{\ell'}^+} \enorm{\rho_{\ell',z}}^2.
\end{equation*}
This concludes the proof.
\end{proof}

\begin{remark}[choice of step-sizes]
By Assumption~\ref{assumption:one}\eqref{assumption:one:a}, the step-sizes need to be uniformly bounded. 
The lower bound for the optimal step-sizes $s_{\ell'}$ follows directly from the finite patch overlap~\eqref{eq:patch_overlap} together with the
definition~\eqref{eq:optimal_step_size}. 
So far the question of a uniform upper bound for $s_{\ell'}$ is still open.
This is precisely why $\lambda_{\ell'}$ is
introduced via definition~\eqref{eq:step_sizes}.
The cut-off bound $d+1$ is arbitrary and could just as well be chosen differently.
The fixed step-size $1/(2(d+1))$, however, is not. 
Its choice is dictated by the finite patch overlap and by Assumption~\ref{assumption:one}\eqref{assumption:one:b}.
Lastly, in our numerical experiments no cut-off takes place and the optimal step-size is always chosen.
\end{remark}

The following theorem presents the main result on the uniform contraction of the geometric multigrid solver for the fractional Laplacian.

\begin{theorem}[uniform contraction of the multigrid for the fractional Laplacian]\label{theorem:main}
Let $p \in \N$ be arbitrary for $d=2$ and $p=1$ otherwise.
Let $u_{\ell}^\star \in \XX_{\ell}^p$ be the discrete solution of~\eqref{eq:discrete_problem} and let $\Phi_{\ell}$ denote the iteration operator of the multigrid solver. 
Then, there exists a factor $0 < \qmg < 1$ such that
\begin{equation}\label{eq:uniform_contraction}
  \enorm{u_\ell^\star - \Phi_{\ell}(u_\ell)} \leq \qmg \, \enorm{u_\ell^\star- u_\ell} \qquad \text{for all } u_{\ell} \in \XX_{\ell}^p.
\end{equation}
The factor $\qmg$ depends only on $s$, $d$, $\TT_0$, $\Omega$, and $\gamma$-shape regularity. 
\end{theorem}

\begin{proof}
  We need to verify Assumptions~\ref{assumption:one}--\ref{assumption:three}.
  Assumption~\ref{assumption:one} follows directly from the statement of Lemma~\ref{lemma:step_size} and $\enorm{\lambda_0 \rho_0}^2 = \lambda_0 \enorm{\rho_0}^2$.
  Assumption~\ref{assumption:three} follows from Proposition~\ref{prop:stable_decomp} for $d \neq 2$ and from Corollary~\ref{cor:hp_stable_decomp} for $d=2$, both of which are shown in section~\ref{sec:stable_decomp}.

  Lastly, we need to show Assumption~\ref{assumption:two}.
  Let $u_0 \in \XX_0^1$, $v_{\ell',z},u_{\ell',z} \in \XX_{\ell',z}^1$ for $\ell' \in \set{1, \dots, \ell-1}$ and $z \in \VV_{\ell'}^+$ as well as $v_{\ell,z} \in \XX_{\ell,z}^p$ for $z \in \VV_{\ell}$ if $p > 1$ and $z \in \VV_{\ell}^+$ if $p=1$.
  We will only treat the first case $p > 1$ in the following.
  Moreover, we set $v_{\ell'} \coloneqq \sum_{z \in \VV_{\ell'}^+}v_{\ell',z}$, $u_k \coloneqq \sum_{w \in \VV_{k}^+} u_{k,w}$, and $v_\ell \coloneqq \sum_{z \in \VV_\ell} v_{\ell,z}$.
  Then, there holds
  \begin{equation}\label{eq:level_decomp}
    \sum_{\ell' = 0}^\ell \sum_{k=0}^{\ell'-1} a(v_{\ell'}, u_k) 
    =
    \sum_{\ell'=1}^{\ell-1}\sum_{k=1}^{\ell'-1} a(v_{\ell'},u_k) + \sum_{\ell'=1}^{\ell-1} a(v_{\ell'},u_0) + \sum_{k=1}^{\ell-1} a(v_\ell,u_k) + a(v_\ell,u_0)
  \end{equation}
  The first term can be estimated using Proposition~\ref{prop:SCS} in section~\ref{sec:SCS}.
  For the last term, the Cauchy--Schwarz inequality and finite patch overlap~\eqref{eq:patch_overlap} yield
  \begin{equation*}
    a(v_\ell,u_0) 
    \leq 
    \enorm{v_\ell} \enorm{u_0} 
    \leq
    \sqrt{2(d+1)} \Big(\sum_{z \in \VV_\ell} \enorm{v_{\ell,z}}^2\Big)^{1/2} \enorm{u_0}.
  \end{equation*}
  With the strengthened Cauchy--Schwarz inequality from Proposition~\ref{prop:SCS} and finite patch overlap, it follows that
  \begin{equation*}
    \Big\lvert \hspace{-0.2em} \Big\lvert \hspace{-0.2em} \Big\lvert \sum_{k=1}^{\ell-1}u_k \Big\rvert \hspace{-0.2em} \Big\rvert \hspace{-0.2em} \Big\rvert^2
    = 
    \sum_{k=1}^{\ell-1} \enorm{u_{k}}^2 + 2 \sum_{\ell'=1}^{\ell-1}\sum_{k=1}^{\ell'-1} a(u_{\ell'},u_k)
    \leq
    2(d+1+\Cscsl)
    \Big(\sum_{k=1}^{\ell-1} \sum_{w \in \VV_{k}^+}\enorm{u_{k,w}}^2\Big)
  \end{equation*}
  and thus
  \begin{align*}
    \sum_{k=1}^{\ell-1} a(v_\ell,u_k) 
    &= 
    a\Big(v_\ell, \sum_{k=1}^{\ell-1}u_k\Big)
    \leq
    \Big\lvert \hspace{-0.2em} \Big\lvert \hspace{-0.2em} \Big\lvert \sum_{k=1}^{\ell-1}u_k \Big\rvert \hspace{-0.2em} \Big\rvert \hspace{-0.2em} \Big\rvert \enorm{v_{\ell}} \\
    &\eqreff*{eq:patch_overlap}{\leq}\,
    2\sqrt{(d+1)(d+1+\Cscsl)} \Big(\sum_{z \in \VV_\ell} \enorm{v_{\ell,z}}^2\Big)^{1/2} \!\Big(\sum_{k=1}^{\ell-1} \sum_{w \in \VV_{k}^+}\enorm{u_{k,w}}^2\Big)^{1/2}.
  \end{align*}
  The second term in~\eqref{eq:level_decomp} can be treated analogously, which implies
  \begin{equation*}
    \sum_{\ell'=0}^{\ell} \sum_{k=0}^{\ell'-1} a(v_{\ell'},u_k) 
    \leq 
    \Cscs
    \Big(\sum_{\ell'=1}^{\ell-1}\sum_{z \in \VV_{\ell'}^+} \enorm{v_{\ell',z}}^2 + \sum_{z \in \VV_{\ell}} \enorm{v_{\ell,z}}^2\Big)^{1/2} \Big(\enorm{u_0}^2+\sum_{k=1}^{\ell-1}\sum_{w \in \VV_{k}^+} \enorm{u_{k,w}}^2\Big)^{1/2}
  \end{equation*}
  with $\Cscs \coloneqq \Cscsl + \sqrt{2(d+1+\Cscsl)}+ 2\sqrt{(d+1)(d+1+\Cscsl)}+\sqrt{2(d+1)}$.
  This concludes the proof of Assumption~\ref{assumption:two}.

  Hence, Theorem~\ref{theorem:uniform_contraction} concludes the proof of uniform contraction~\eqref{eq:uniform_contraction}.
\end{proof}

\subsection{Multigrid as inexact solver for AFEM}

As in~\cite{IMPS2024}, the multigrid solver from Algorithm~\ref{alg:multigrid} can be employed as an algebraic solver within an adaptive algorithm, since Theorem~\ref{theorem:main} guarantees that it is uniformly contractive.
As a direct consequence of this, in particular, the computable consecutive iterate difference $\enorm{u_\ell^k-u_\ell^{k-1}}$ provides a two-sided estimate of the algebraic error and can thus be used in the stopping criterion.
Combined with nested iteration and an adaptive stopping criterion, the resulting Algorithm~\ref{alg:afem} below provably achieves optimal convergence rates; see~\cite{ghps2021}.
Its remaining modules are standard: elements are selected by the Dörfler marking criterion~\cite{doerfler1996}, and meshes are refined by NVB.
To estimate the discretization error, we use the weighted residual error estimator from~\cite{FMP2021} introduced for lowest-order elements, i.e., $p=1$, which satisfies the axioms of adaptivity from~\cite{axioms}.

\begin{definition}[weighted residual error estimator \cite{FMP2021}]\label{def:estimator}
Consider a mesh $\mathcal T_\ell$.
For $f\in L^2(\Omega)$ and $v_\ell\in \XX_{\ell}^1$, extended by zero to $\mathbb R^d$, we define the refinement indicators
\begin{equation}\label{eq:estimator:local}
\eta_\ell(T;v_\ell) :=
\begin{cases}
\big\|h_T^{s}\,\big(f-(-\Delta)^s v_\ell\big)\big\|_{L^2(T)} & \text{for } 0<s\le 1/2,\\[4pt]
\big\|h_T^{1/2}\operatorname{dist}(\cdot,\partial T)^{s-1/2}\,\big(f-(-\Delta)^s v_\ell\big)\big\|_{L^2(T)} & \text{for } 1/2<s<1,
\end{cases}
\end{equation}
for all $T\in\mathcal T_\ell$, where $h_T:=|T|^{1/d}$. Moreover, for all $\mathcal U_\ell\subseteq\mathcal T_\ell$, we set
\begin{equation}\label{eq:estimator:global}
\eta_\ell(\mathcal U_\ell;v_\ell) := \Big(\sum_{T\in\mathcal U_\ell}\eta_\ell(T;v_\ell)^2\Big)^{1/2}
\quad\text{and}\quad
\eta_\ell(v_\ell) := \eta_\ell(\mathcal T_\ell;v_\ell).
\end{equation}
\end{definition}

 The adaptive algorithm with inexact solver is given as follows.
 
\begin{algorithm}[AFEM with iterative solver]\label{alg:afem}
\textbf{Input:} initial mesh $\mathcal T_0$, 
parameters $0<\theta\le1$, $C_{\mathrm{mark}}\ge1$, $\mu>0$, initial guess $u_0^0:=0$.\\
\textbf{Adaptive loop:} For all $\ell=0,1,2,\dots$, repeat the following steps (I)--(III):
\begin{itemize}
\item[(I)] \textup{\texttt{SOLVE \& ESTIMATE:}} For all $k=1,2,3,\dots$, repeat the following steps (i)--(iii):
  \begin{itemize}
  \item[(i)] Compute $u_\ell^k := \Phi(u_\ell^{k-1})$ via Algorithm~\ref{alg:multigrid}.
  \item[(ii)] Compute the refinement indicators $\{\eta_\ell(T;u_\ell^k)\}_{T\in\mathcal T_\ell}$ as in Definition~\ref{def:estimator}.
  \item[(iii)] If $\enorm{u_\ell^k-u_\ell^{k-1}}\le\mu\,\eta_\ell(u_\ell^k)$, terminate the $k$-loop and set $\underline k[\ell]:=k$ and $u_\ell:=u_\ell^{\underline k[\ell]}$.
  \end{itemize}
\item[(II)] \textup{\texttt{MARK:}} Determine a set $\mathcal M_\ell\subseteq\mathcal T_\ell$ of up to the factor $C_{\mathrm{mark}}$ minimal cardinality such that $\theta\,\eta_\ell(u_\ell)^2\le\eta_\ell(\mathcal M_\ell;u_\ell)^2$.
\item[(III)] \textup{\texttt{REFINE:}} Generate $\mathcal T_{\ell+1}:=\textup{\texttt{REFINE}}(\mathcal T_\ell,\mathcal M_\ell)$ by NVB and set $u_{\ell+1}^0:=u_\ell$.
\end{itemize}
\textbf{Output:} meshes $\mathcal T_\ell$, approximations $u_\ell$, and estimators $\eta_\ell(u_\ell)$.
\end{algorithm}

\begin{remark}[computational cost] Provided that each step of the adaptive Algorithm~\ref{alg:afem} can be implemented in linear complexity, and since the multigrid solver is shown to be uniformly contractive in Theorem~\ref{theorem:main}, then \cite{ghps2021} ensures optimal complexity, i.e., optimal convergence rates with respect to the total computational cost and, hence, time. For the nonlocal operator at hand, however, the Galerkin matrices are dense, so that neither a solver step nor the evaluation of the estimator is currently performed in linear complexity. More precisely, while the structure of the solver remains the same as in~\cite{IMPS2024}, the computation of the residual in~\eqref{eq:def_local_updates} requires a matrix-vector product with a dense matrix. When optimal step-sizes are used, their computation involves a further matrix-vector product of the same type. These dense matrix-vector products are the main computational bottleneck for the solver. Overcoming it requires matrix compression techniques such as hierarchical matrices, which is beyond the scope of this work.
Nonetheless, the numerical experiments in Section~\ref{sec:num} show that employing the robust solver within the adaptive algorithm already reduces the computational effort substantially.
\end{remark}

For $r>0$, the approximation class is characterized by
\begin{equation*}
\|u^\star\|_{\mathbb A_r}:=\sup_{N\in\mathbb N_0}\,(N+1)^r\min_{\mathcal T_{\mathrm{opt}}\in\mathbb T_N}\big(\enorm{u^\star-u_{\mathrm{opt}}^\star}+\eta_{\mathrm{opt}}(u_{\mathrm{opt}}^\star)\big),
\end{equation*}
where $\mathbb T_N:=\{\mathcal T\in\mathbb T:\#\mathcal T-\#\mathcal T_0\le N\}$ and $u_{\mathrm{opt}}^\star$ 
is the piecewise linear Galerkin solution on $\mathcal T_{\mathrm{opt}}$.
Thus, $\|u^\star\|_{\mathbb A_r}<\infty$ means that, on optimally chosen meshes, the error plus estimator decays with rate $r$ with respect to the number of elements (and thus number of degrees of freedom).
The following theorem shows that Algorithm~\ref{alg:afem} achieves precisely these rates.
 
\begin{theorem}[optimal convergence rates]\label{thm:afem}
Consider the index set $\mathcal Q:=\{(\ell,k)\in\mathbb N_0^2 : u_\ell^k \text{ is computed in Algorithm~\ref{alg:afem}}\}$ and, for all $(\ell,k)\in\mathcal Q$, define the quasi-error
\begin{equation*}
\Delta_\ell^k:=\enorm{u^\star-u_\ell^k}+\eta_\ell(u_\ell^k).
\end{equation*}
There exist $0<\theta^\star\le1$ and $\mu^\star>0$ such that, for all $0<\theta<\theta^\star$, $0<\mu/\theta<\mu^\star$, and $r>0$, there holds
\begin{equation*}
c_{\mathrm{opt}}\,\|u^\star\|_{\mathbb A_r}
\le\sup_{(\ell,k)\in\mathcal Q}(\#\mathcal T_\ell)^r\,\Delta_\ell^k
\le C_{\mathrm{opt}}\max\big\{\|u^\star\|_{\mathbb A_r},\Delta_0^0\big\}.
\end{equation*}
The constants $c_{\mathrm{opt}},C_{\mathrm{opt}}>0$ depend only on the initial triangulation $\mathcal T_0$, the fractional order $s$, the rate $r$, the marking constant $C_{\mathrm{mark}}$, the ratios $\theta/\theta^\star$ and $\mu/(\theta\mu^\star)$, and the properties of NVB.
\end{theorem}
 
\begin{proof}
The estimator from Definition~\ref{def:estimator} satisfies the axioms of adaptivity~\cite{axioms}; see~\cite{FMP2021}. Together with the uniform contraction and the two-sided algebraic error estimate from Theorem~\ref{theorem:main}, the result follows from~\cite[Thm.~8]{ghps2021}; see also~\cite[Sect.~3]{IMPS2024}.
\end{proof}
 

\section{Strengthened Cauchy--Schwarz inequality}\label{sec:SCS}

To show the required strengthened Cauchy--Schwarz inequality for $\widetilde{H}^s(\Omega)$, we rely both on the ideas from~\cite{IMPS2024} to prove said inequality for $H_0^1(\Omega)$ as well as a splitting into near- and far-field interactions to deal with the non-local energy norm. 
We first introduce some auxiliary results.

\subsection{Auxiliary results}
For every level $\ell \in \N_0$ and $z \in \VV_\ell$, define the generation of the patch $\TT_{\ell}(z)$ by
\begin{equation*}
  g_{\ell,z} \coloneqq \max_{T \in \TT_{\ell}(z)} \operatorname{level}(T) 
  \coloneqq \max_{T \in \TT_{\ell}(z)} \log_2(\abs{T_0}/\abs{T}) \in \N_0,
\end{equation*}
where $T_0 \in \TT_0$ is the unique ancestor of $T \in \TT_{\ell}(z)$.
Furthermore, denote the smallest level on the two-patch via
\begin{equation}\label{def:smallest_level}
  r_{\ell,z} \coloneqq \min_{T \in \TT_{\ell}^2(\{z\})} \operatorname{level}(T).
\end{equation}
There holds the following lemma for the generation and $r_{\ell,z}$. 
We refer to~\cite[Lemma~10 \& Lemma~9]{FFPS2017} and~\cite[Lemma 22]{hmp26} for a proof.
\begin{lemma}\label{lemma:uniform_patch}
 There exists a constant $\Cspan \in \N$ depending only on $\gamma$-shape regularity, such that
\begin{equation}\label{eq:level_span}
  g_{\ell,z} \leq r_{\ell-1,z} + \Cspan \quad \text{for all } \ell \in \N.
\end{equation}
Moreover, there exists a uniform $n \in \N_0$ depending only on the initial mesh $\TT_0$ such that 
\begin{equation}\label{eq:uniform_patch}
  \overline{\omega_{\ell}(z)} \subseteq \omega_{\ell-1}^2(\{z\}) \subseteq \widehat{\omega}_{g_{\ell,z}}^n(z) \qquad \text{for all } \ell\in \N \text{ and all } z \in \VV_\ell.
\end{equation}
\end{lemma}

Next, we show a Poincaré inequality for $\widetilde{H}^s(\Omega)$. A fractional Poincaré inequality can, e.g., be found in \cite{AB2017}; for sake of completeness and in order to make the domain scaling explicit, we provide a short proof. 

\begin{lemma}[Poincaré inequality]\label{lemma:local_estimate}
  Let $v_{\ell} \in \XX_{\ell}^p$ and $\omega \subset \Omega$ be a bounded subset with $\supp v_{\ell} \subseteq \omega$ and $h \coloneqq \diam \omega \coloneqq \sup_{x,y\in \omega} |x-y|$.
  Then, there holds 
  \begin{equation}\label{eq:local_estimate}
    \norm{v_{\ell}}_{L^2(\Omega)} \leq \Cloc \, h^s \, \enorm{v_{\ell}}.
  \end{equation}
  The constant $\Cloc$ depends only on $d$ and $s$.
\end{lemma}

\begin{proof}
Let $x \in \omega$ and define $A_x \coloneqq \set{y \in \R^d \given 2h < \abs{x-y} < 3h}$. 
For $y \in A_x$ and any $z \in \omega$, it follows that $\abs{z-y} \geq \abs{x-y} - \abs{z-x} \geq 2h-h \geq h$. 
This implies $y \notin \omega$ and hence $v_{\ell}(y) = 0$.
Therefore, we obtain
\begin{equation*}
  \abs{v_{\ell}(x)}^2 = \frac{1}{\abs{A_x}} \int_{A_x} \abs{v_\ell(x)-v_\ell(y)}^2 \d{y} 
  \leq 
  \frac{(3h)^{d+2s}}{\abs{A_x}} \int_{\R^d} \frac{\abs{v_\ell(x)-v_\ell(y)}^2}{\abs{x-y}^{d+2s}} \d{y}.
\end{equation*}
Since $\abs{A_x} = (3^d-2^d) \, h^d \, \abs{\S_d}$, where $\abs{\S_d}$ denotes the volume of the unit sphere in $\R^d$, we conclude 
\begin{equation*}
  \norm{v_{\ell}}^2_{L^2(\Omega)} \leq \frac{(3h)^{d+2s}}{(3^d-2^d) h^d \abs{\S_d}} \int_{\R^d}\int_{\R^d} \frac{\abs{v_\ell(x)-v_\ell(y)}^2}{\abs{x-y}^{d+2s}} \d{x} \d{y} = \frac{3^{d+2s}2}{(3^d-2^d) \abs{\S_d} C(d,s)} \, h^{2s} \, \enorm{v_\ell}^2
\end{equation*} 
and set $\Cloc^2 \coloneqq 3^{d+2s}2/\big((3^d-2^d)\,\abs{\S_d}\,C(d,s)\big)$.
\end{proof}

\subsection{Proof of strengthened Cauchy--Schwarz inequality}

Since we use generation constraints in the proof, the following notation will be helpful. 
For $i \in \N_0$, we define the set of new vertices and their direct neighbors with a fixed generation via $\VV_{\ell}^+(i) \coloneqq \set{z \in \VV_{\ell}^+ \given g_{\ell,z} = i}.$
With this we can state and prove the main result of this section.

\begin{proposition}[strengthened Cauchy--Schwarz inequality]\label{prop:SCS}
Let $\ell \in \N_0$. 
Let $v_{\ell',z} \in \XX_{\ell',z}^1$ for $\ell' = 1, \dots, \ell-1$ and $z \in \VV_{\ell'}^+$ as well as $u_{k,w} \in \XX_{k,w}^1$ for $k = 1, \dots, \ell-2$ and $w \in \VV_k^+$. 
Then, there exists a constant $\Cscsl > 0$ such that
\begin{equation}\label{eq:SCS}
 \sum_{\ell'=1}^{\ell-1}\sum_{k=1}^{\ell'-1} a\Big(\sum_{z \in \VV_{\ell'}^+} v_{\ell',z}, \sum_{w \in \VV_k^+} u_{k,w}\Big)
  \leq 
  \Cscsl \, \Big(\sum_{\ell'=1}^{\ell-1}\sum_{z \in \VV_{\ell'}^+} \enorm{v_{\ell',z}}^2\Big)^{1/2}
  \Big(\sum_{k=1}^{\ell-2}\sum_{w \in \VV_k^+} \enorm{u_{k,w}}^2\Big)^{1/2}.
\end{equation}
The constant \( \Cscsl > 0\) depends only on the initial mesh $\TT_0$, $\Omega$, $d$, $s$ and $\gamma$-shape regularity.
\end{proposition}

\begin{proof}
This proof consists of four steps, following the strategy outlined below.
First, we split the functions $v_{\ell',z}$ and $u_{k,w}$ into two groups: those whose supports lie far apart, which are called \emph{far-field interactions}, and those whose supports lie close together, which are called \emph{near-field interactions}. 
Provided ``far apart'' is defined appropriately, the far-field interactions can be treated directly via the definition of the bilinear form $a(\cdot, \cdot)$. 
The near-field interactions, in turn, require additional care due to the non-local norm, but can otherwise be handled similarly to the $H_0^1(\Omega)$ case; see, e.g., the proof of~\cite[Proposition 5.7]{IMPS2024}.

\textbf{Step 1 (decomposing the sums):}
Let $M \coloneqq \max_{z \in \VV_\ell} g_{\ell,z}$.
For any $i,j \in \{0, \dots, M\}$, $z \in \widehat{\VV}_i$, and $w \in \widehat{\VV}_j$, and with $n \in \N_0$ given by Lemma~\ref{lemma:uniform_patch}, define the interaction distance by
\begin{equation}\label{def:distance}
  D_{ij}^{zw}
  \coloneqq \dist{\widehat{\omega}^n_i(z)}{\widehat{\omega}^n_j(w)} \coloneqq \inf_{\substack{x \in \widehat{\omega}^n_i(z) \\ y \in \widehat{\omega}^n_j(w)}} \abs{x-y}.
\end{equation}
Note that $D_{ij}^{zw} = D_{ji}^{wz}$ and set $h_{ij} \coloneqq \max\set{\widehat{h}_i,\widehat{h}_j}$.
Let $\ell' \in \set{1, \dots, \ell-1}$, $z \in \VV_{\ell'}^+(i)$, $k \in \set{1,\dots,\ell-2}$, and $w \in \VV_k^+(j)$.
We refer to the pairs $(\ell', z)$ and $(k, w)$ as a \emph{far-field interaction}, if the admissibility condition
\begin{equation}
D_{ij}^{zw} > h_{ij}
\end{equation}
holds.
Otherwise, we call them \emph{near-field interactions}.
Rewriting the left-hand side of~\eqref{eq:SCS} by introducing generation constraints results in
\begin{equation}\label{eq:SCS:generation_constraints}
  \sum_{\ell'=1}^{\ell-1}\sum_{k=1}^{\ell'-1}\sum_{z \in \VV_{\ell'}^+}\sum_{w \in \VV_k^+} a(v_{\ell',z}, u_{k,w}) 
  =
  \sum_{i=0}^M \sum_{j=0}^M \sum_{\ell'=1}^{\ell-1}\sum_{k=1}^{\ell'-1}\sum_{z \in \VV_{\ell'}^+(i)}\sum_{w \in \VV_k^+(j)} a(v_{\ell',z}, u_{k,w}).
\end{equation}
We decompose~\eqref{eq:SCS:generation_constraints} into far-field and near-field interactions and split the near-field interactions into the two cases $i \leq j$ and $j < i$. 
In the latter case, exchanging the roles of the summation pairs $(\ell', z)$ and $(k, w)$ yields the identity
\begin{equation*} 
  \sum_{i=0}^M \sum_{j=0}^{i-1} \sum_{\ell'=1}^{\ell-1}\sum_{z \in \VV_{\ell'}^+(i)} \sum_{k=1}^{\ell'-1} \sum_{\substack{w \in \VV_k^+(j) \\ D_{ij}^{zw} \leq h_{ij}}} a(v_{\ell',z}, u_{k,w})  
  =  
  \sum_{j=0}^M \sum_{i=j+1}^{M} \sum_{k=1}^{\ell-2}\sum_{w \in \VV_{k}^+(j)} \sum_{\ell'=k+1}^{\ell-1} \sum_{\substack{z \in \VV_{\ell'}^+(i) \\ D_{ji}^{wz} \leq h_{ij}}} a(v_{\ell',z}, u_{k,w}), 
\end{equation*}
Substituting this identity back into~\eqref{eq:SCS:generation_constraints} gives
\begin{align*}
  &\sum_{i=0}^M \sum_{j=0}^M \sum_{\ell'=1}^{\ell-1}\sum_{k=1}^{\ell'-1}\sum_{z \in \VV_{\ell'}^+(i)} \sum_{w \in \VV_k^+(j)} \!\! a(v_{\ell',z}, u_{k,w})  
  = 
  \sum_{i=0}^M \sum_{j=0}^M \sum_{\ell'=1}^{\ell-1}\sum_{z \in \VV_{\ell'}^+(i)} \sum_{k=1}^{\ell'-1} \sum_{\substack{w \in \VV_k^+(j) \\ D_{ij}^{zw} > h_{ij}}} \!\! a(v_{\ell',z}, u_{k,w}) \\ 
  & \; + 
  \sum_{i=0}^M \sum_{j=i}^M \sum_{\ell'=1}^{\ell-1}\sum_{z \in \VV_{\ell'}^+(i)} \sum_{k=1}^{\ell'-1} \sum_{\substack{w \in \VV_k^+(j) \\ D_{ij}^{zw} \leq h_{ij}}} \!\! a(v_{\ell',z}, u_{k,w})  
  +  
  \sum_{j=0}^M \sum_{i=j+1}^{M} \sum_{k=1}^{\ell-2}\sum_{w \in \VV_{k}^+(j)} \sum_{\ell'=k+1}^{\ell-1} \sum_{\substack{z \in \VV_{\ell'}^+(i) \\ D_{ji}^{wz} \leq h_{ij}}} \!\! a(v_{\ell',z}, u_{k,w}) \\
  & \qquad \qquad \eqqcolon 
  \sum_{i=0}^M\sum_{j=0}^M S_1(i,j) + \sum_{i=0}^M\sum_{j=i}^{M} S_2(i,j) + \sum_{j=0}^{M}\sum_{i=j+1}^M S_3(j,i).
\end{align*}

\textbf{Step 2 (far-field interactions):}
Let the generations $i,j \in \{0, \dots, M\}$ be fixed and let us consider the sum $S_1(i,j)$. 
Recall $\supp v_{\ell',z} \subseteq \widehat{\omega}_i^n(z)$ and $\supp u_{k,w} \subseteq \widehat{\omega}_j^n(w)$, and since $\widehat{\omega}_i^n(z) \cap \widehat{\omega}_j^n(w) = \emptyset$ holds for far-field interactions, we have
\begin{equation*}
  a(v_{\ell',z}, u_{k,w}) 
  = 
  - C(d,s) \int_{\widehat{\omega}_i^n(z)} \int_{\widehat{\omega}_j^n(w)} \!\! \frac{v_{\ell',z}(x)u_{k,w}(y)}{\abs{x-y}^{d+2s}} \d{x}\d{y} 
  \eqreff{def:distance}\leq 
  \frac{C(d,s)}{(D_{ij}^{zw})^{d+2s}} \, \norm{v_{\ell',z}}_{L^1(\Omega)} \, \norm{u_{k,w}}_{L^1(\Omega)}.
\end{equation*}
Using the Cauchy--Schwarz inequality and Lemma~\ref{lemma:local_estimate}, it follows that
\begin{equation*}
  \norm{v_{\ell',z}}_{L^1(\Omega)} 
  \leq 
  \abs{\widehat{\omega}_i^n(z)}^{1/2} \norm{v_{\ell',z}}_{L^2(\Omega)}
  \eqreff{eq:local_estimate}{\lesssim} 
  \,\widehat{h}_i^{d/2+s} \enorm{v_{\ell',z}}
\end{equation*}
and hence also
\begin{equation}\label{eq:SCS:1}
  a(v_{\ell',z},u_{k,w}) 
  \lesssim 
  \frac{(\widehat{h}_i\widehat{h}_j)^{(d+2s)/2}}{(D_{ij}^{zw})^{d+2s}} \, \enorm{v_{\ell',z}} \, \enorm{u_{k,w}}.
\end{equation}
Applying this estimate and the discrete Cauchy--Schwarz inequality in the form
\begin{equation*}
  \sum_{\alpha \in I} x_{\alpha}y_{\alpha} 
  \leq 
  \Big(\sum_{\alpha \in I} x_{\alpha}^2\Big)^{1/2} \Big(\sum_{\alpha \in I} y_{\alpha}^2\Big)^{1/2}
\end{equation*}
to $S_1(i,j)$ with $I \coloneqq \set{(\ell',z,k,w) \given 1 \leq \ell' <\ell, z \in \VV_{\ell'}^+(i), 1 \leq k < \ell' , w \in \VV_k^+(j), D_{ij}^{zw} > h_{ij}}$ yields
\begin{align}\label{eq:SCS:discrete_CS}
  \begin{split}
  \sum_{\ell'=1}^{\ell-1}\sum_{z \in \VV_{\ell'}^+(i)} &\sum_{k=1}^{\ell'-1}\sum_{\substack{w \in \VV_{k}^+(j) \\ D_{ij}^{zw}>h_{ij}}} a(v_{\ell',z}, u_{k,w}) 
  \eqreff{eq:SCS:1}{\lesssim} 
  \sum_{\ell'=1}^{\ell-1}\sum_{z \in \VV_{\ell'}^+(i)} \sum_{k=1}^{\ell'-1}\sum_{\substack{w \in \VV_{k}^+(j) \\ D_{ij}^{zw}>h_{ij}}} \frac{(\widehat{h}_i\widehat{h}_j)^{(d+2s)/2}}{(D_{ij}^{zw})^{d+2s}} \, \enorm{v_{\ell',z}} \, \enorm{u_{k,w}} \\
  &\leq 
  (\widehat{h}_i\widehat{h}_j)^{(d+2s)/2}
  \Big(\sum_{\ell'=1}^{\ell-1} \sum_{z \in \VV_{\ell'}^+(i)} \sum_{k=1}^{\ell'-1}\sum_{\substack{w \in \VV_{k}^+(j) \\ D_{ij}^{zw}>h_{ij}}}
  (D_{ij}^{zw})^{-d-2s} \, \enorm{v_{\ell',z}}^2 \Big)^{1/2}  \\
  &\hspace{4cm} \Big(\sum_{\ell'=1}^{\ell-1}  \sum_{z \in \VV_{\ell'}^+(i)} \sum_{k=1}^{\ell'-1}\sum_{\substack{w \in \VV_{k}^+(j) \\ D_{ij}^{zw}>h_{ij}}} (D_{ij}^{zw})^{-d-2s} \, \enorm{u_{k,w}}^2 \Big)^{1/2}.
  \end{split}
\end{align}
By further bounding both factors using $\ell' \leq \ell-1$, and by exchanging the summation order of the pairs $(\ell',z)$ and $(k,w)$ in the second factor, we obtain
\begin{align}\label{eq:SCS:separate_sums}
  \sum_{\ell'=1}^{\ell-1} \sum_{z \in \VV_{\ell'}^+(i)} \sum_{k=1}^{\ell'-1}\sum_{\substack{w \in \VV_{k}^+(j) \\ D_{ij}^{zw}>h_{ij}}}  (D_{ij}^{zw})^{-d-2s} \, \enorm{v_{\ell',z}}^2
  &\leq 
  \sum_{\ell'=1}^{\ell-1} \sum_{z \in \VV_{\ell'}^+(i)} \enorm{v_{\ell',z}}^2 
  \sum_{k=1}^{\ell-2}\sum_{\substack{w \in \VV_{k}^+(j) \\ D_{ij}^{zw}>h_{ij}}} (D_{ij}^{zw})^{-d-2s}, \nonumber \\
  \sum_{\ell'=1}^{\ell-1} \sum_{z \in \VV_{\ell'}^+(i)}  \sum_{k=1}^{\ell'-1}\sum_{\substack{w \in \VV_{k}^+(j) \\ D_{ij}^{zw}>h_{ij}}}  (D_{ij}^{zw})^{-d-2s} \, \enorm{u_{k,w}}^2
  &\leq 
  \sum_{k=1}^{\ell-2} \sum_{w \in \VV_{k}^+(j)} \enorm{u_{k,w}}^2 \sum_{\ell'=1}^{\ell-1}\sum_{\substack{z \in \VV_{\ell'}^+(i) \\ D_{ji}^{wz}>h_{ij}}} (D_{ij}^{zw})^{-d-2s}.
\end{align}
We only treat the first factor and the second one follows analogously.
In~\cite[Lemma 3.1]{wc2006} it is shown, for $d=2$, that there exists a constant $\Clev$ depending only on $\gamma$-shape regularity, such that
\begin{equation}\label{eq:SCS:level_bound}
  \# \LL_{\underline{\ell}}^{\overline{\ell}}(w,j) 
  \coloneqq 
  \# \set{k \in \set{\underline{\ell}, \dots, \overline{\ell}} \given  w \in \VV_{k}^+(j)} 
  \leq 
  \Clev < \infty \qquad \text{for all } w \in \VV_{\ell}.
\end{equation}
Since the arguments in~\cite[Lemma 3.1]{wc2006} are based on $\gamma$-shape regularity, this result can be extended to any dimension.  
Moreover, for far-field interactions, i.e., in the case where $D_{ij}^{zw} > h_{ij}$, there holds
\begin{equation}\label{eq:SCS:distance_equivalence}
  \abs{z-w}
  \eqreff{def:distance}{\leq} 
  D_{ij}^{zw} + \diam(\widehat{\omega}^n_i(z))+ \diam(\widehat{\omega}^n_j(w))
  \simeq 
  D_{ij}^{zw} + \widehat{h}_i + \widehat{h}_j 
  \lesssim 
  D_{ij}^{zw}
  \eqreff{def:distance}{\leq} 
  \abs{z-w}.
\end{equation}
Recall that $i,j \in \{1,\dots, M\}$ are fixed and consider a fixed $z \in \widehat{\VV}_{i}$.
The above observations yield
\begin{align*}
  \sum_{k=1}^{\ell-2}\sum_{\substack{w \in \VV_{k}^+(j) \\ D_{ij}^{zw}>h_{ij}}}  (D_{ij}^{zw})^{-d-2s} 
  &=  
  \sum_{\substack{w \in \widehat{\VV}_j \\ D_{ij}^{zw} > h_{ij}}} \sum_{k \in \LL_1^{\ell-2}(w,j)}  (D_{ij}^{zw})^{-d-2s} \\
  &\eqreff*{eq:SCS:level_bound}{\leq} \;
  \Clev \sum_{\substack{w \in \widehat{\VV}_j \\ D_{ij}^{zw} > h_{ij}}}  (D_{ij}^{zw})^{-d-2s} 
  \eqreff{eq:SCS:distance_equivalence}{\lesssim}
  \sum_{\substack{w \in \widehat{\VV}_j \\ \abs{z-w} > h_{ij}}} \abs{z-w}^{-d-2s}. 
\end{align*}
Next, we want to establish a uniform bound for this sum. 
As the mesh $\widehat{\TT}_j$ is quasi-uniform, there holds $\abs{w-w'} \gtrsim \widehat{h}_j$ for all $w \neq w' \in \widehat{\VV}_j$.
Hence, the number of vertices $w \in \widehat{\VV}_j$ contained in any ball of radius $r$ around $z$ can be bounded by $C \, (r/\widehat{h}_j)^d$, where $C$ depends only on the initial mesh $\TT_0$ and the dimension $d$.
Summing over $w \in \widehat{\VV}_j$ with $|z-w| > h_{ij}$ is equivalent to summing over $n \in \N_0$ and the annuli $A_n^{z} \coloneqq \set{w \in \widehat{\VV}_j \given 2^nh_{ij} < \abs{z-w} \leq 2^{n+1} h_{ij}}$ with $\abs{A_n^z} \lesssim (2^{n+1}h_{ij}/\widehat{h}_j)^d$.
Altogether, we obtain
\begin{align}\label{eq:SCS:annuli}
  \sum_{\substack{w \in \widehat{\VV}_j \\ \abs{z-w} > h_{ij}}} \!\!\! &\abs{z-w}^{-d-2s}
  = 
  \sum_{n \in \N_0} \sum_{w \in A_n^{z}} \abs{z-w}^{-d-2s}
  \leq
  \sum_{n \in \N_0} \sum_{w \in A_n^{z}} (2^n h_{ij})^{-d-2s}  \leq
  \sum_{n \in \N_0} \abs{A_n^{z}} (2^n h_{ij})^{-d-2s} \nonumber \\
  &\lesssim
  \sum_{n \in \N_0} \Big(\frac{2^{n+1}h_{ij}}{\widehat{h}_j}\Big)^d (2^n h_{ij})^{-d-2s} 
  = 
  2^d h_{ij}^{-2s}\widehat{h}_j^{-d} \sum_{n \in \N_0} 2^{-2sn}
  \lesssim h_{ij}^{-2s}\widehat{h}_j^{-d},
\end{align}
where the remaining sum is bounded as a geometric series and the hidden constant depends only on $\TT_0$, $\gamma$-shape regularity, $d$, and $s$.
The right-hand side in~\eqref{eq:SCS:separate_sums} can be estimated in the same way.
Thus,~\eqref{eq:SCS:annuli} together with~\eqref{eq:SCS:discrete_CS} establishes
\begin{equation*}
  S_1(i,j)
  \lesssim
  (\widehat{h}_i\widehat{h}_j)^{(d+2s)/2} h_{ij}^{-2s}\widehat{h}_j^{-d/2}\widehat{h}_i^{-d/2} \, \Big(\sum_{\ell'=1}^{\ell-1}\sum_{z \in \VV_{\ell'}^+(i)} \enorm{v_{\ell',z}}^2\Big)^{1/2} 
  \Big(\sum_{k=1}^{\ell-2}\sum_{w \in \VV_k^+(j)} \enorm{u_{k,w}}^2\Big)^{1/2}. 
\end{equation*}
Rewriting the mesh-size functions yields
\begin{align*}
  (\widehat{h}_i\widehat{h}_j)^{(d+2s)/2} h_{ij}^{-2s} \widehat{h}_j^{-d/2}\widehat{h}_i^{-d/2}
  &= 
  (\widehat{h}_i\widehat{h}_j)^{s} h_{ij}^{-2s} 
  =
  (\min\{\widehat{h}_i,\widehat{h}_j\} \max\{\widehat{h}_i,\widehat{h}_j\})^{s} h_{ij}^{-2s}\\
  &=
  \Big(\frac{\min\{\widehat{h}_i,\widehat{h}_j\}}{\max\{\widehat{h}_i,\widehat{h}_j\}}\Big)^s
  = 
  2^{-|i-j|s/d}
\end{align*}
and thus
\begin{equation}\label{eq:SCS:far_field_result}
  S_1(i,j) 
  \lesssim 
  2^{-\abs{i-j}s/d} \Big(\sum_{\ell'=1}^{\ell-1}\sum_{z \in \VV_{\ell'}^+(i)} \enorm{v_{\ell',z}}^2\Big)^{1/2} 
  \Big(\sum_{k=1}^{\ell-2}\sum_{w \in \VV_k^+(j)} \enorm{u_{k,w}}^2\Big)^{1/2}, 
\end{equation}
which concludes this step.

\textbf{Step 3 (near-field interactions):} 
We note that the proofs for the contributions $S_2(i,j)$ and $S_3(j,i)$ follow with identical arguments. Therefore, we only consider $S_2$ in the following.
Let the generations $i,j \in \{1, \dots, M\}$ be fixed and assume $i \leq j$.
In this step, we are treating the following sum $S_2(i,j)$, where 
\begin{equation*}
  S_2(i,j)
  =
  \sum_{\ell'=1}^{\ell-1}\sum_{z \in \VV_{\ell'}^+(i)} \sum_{k=1}^{\ell'-1} \sum_{\substack{w \in \VV_k^+(j) \\ D_{ij}^{zw} \leq h_{ij}}} a(v_{\ell',z}, u_{k,w}) 
  = 
  \sum_{\ell'=1}^{\ell-1}\sum_{z \in \VV_{\ell'}^+(i)} a\Big(v_{\ell',z}, \sum_{k=1}^{\ell'-1} \sum_{\substack{w \in \VV_k^+(j) \\ D_{ij}^{zw} \leq h_{ij}}} u_{k,w}\Big).
\end{equation*}
Due to the generation constraint $w \in \VV_k^+(j)$, there holds $\sum_{k=1}^{\ell'-1} \sum_{w \in \VV_k^+(j), D_{ij}^{zw} \leq h_{ij}} u_{k,w} \in \widehat{\XX}_j^1$.
Because $i \leq j$, we can apply~\cite[Lemma 5.7]{FMP2021_SZ} for quasi-uniform meshes to obtain
\begin{equation}\label{eq:SCS:uniformSCS}
  S_2(i,j) 
  \leq 
  \Cuscs \, \delta^{j-i} \,
  \sum_{\ell'=1}^{\ell-1}\sum_{z \in \VV_{\ell'}^+(i)} \widehat{h}_j^{-s} \,\enorm{v_{\ell',z}} \, \norm[\Big]{\sum_{k=1}^{\ell'-1} \sum_{\substack{w \in \VV_k^+(j) \\ D_{ij}^{zw} \leq h_{ij}}} u_{k,w}}_{L^2(\Omega)},
\end{equation}
where $\delta \coloneqq 2^{-(s-\beta)/d}$ with $\beta \coloneqq \max \{s-1/2,0\}$ and $\Cuscs$ depends only on $\Omega$, $d$, $s$, and the initial mesh $\TT_0$.
Consider first the $L^2(\Omega)$-norm on the right-hand side of the above inequality. 
Let $x \in \Omega$ and $T_x \in \widehat{\TT}_j$ be one element such that $x \in T_x$.
Since $u_{k,w} \in \XX_{k,w}^1$, then $\supp u_{k,w} \subseteq \overline{\omega_k(w)}$ and only the functions $u_{k,w}$ with $x \in \overline{\omega_k(w)}$ can be non-zero at $x$.
If $x \in \overline{\omega_k(w)}$, then at least one of the vertices $(z^x_m)_{m=1}^{d+1}$ of $T_x$ also lies in $\overline{\omega_k(w)}$.
Let 
\begin{equation*}
  \mathcal{LV}(j,z_m^x)
  \coloneqq 
  \set{(k,w) \given k \in \{1,\dots,\ell'-1\}, w \in \VV_k^+(j), \text{ and } z^x_m \in \overline{\omega_k(w)}}.
\end{equation*}
Moreover, in~\cite[Lemma 3.1]{wc2006} it is shown, for $d=2$, that there exists a constant $\Clevv$ depending only on shape-regularity, such that
\begin{equation*}
  \# \mathcal{LV}(j,z_m^x) \leq \Clevv < \infty.
\end{equation*}
As in the previous step, since the arguments in~\cite[Lemma 3.1]{wc2006} are based on $\gamma$-shape regularity, this result can be extended to any dimension. Since 
\begin{align*}
\# \set{(k,w) \given  1 \leq k < \ell', w \in \VV_k^+(j), D_{ij}^{zw} \leq h_{ij},u_{k,w}(x) \neq 0} \leq \sum_{m=1}^{d+1} \# \mathcal{LV}(j,z_m^x) \leq (d+1) \Clevv,
\end{align*}
with the discrete Cauchy--Schwarz inequality, this implies
\begin{align*}
  \abs[\Big]{\sum_{k=1}^{\ell'-1} \sum_{\substack{w \in \VV_k^+(j) \\ D_{ij}^{zw} \leq h_{ij}}} &u_{k,w}(x)}^2 
  \leq 
  (d+1) \, \Clevv \sum_{k=1}^{\ell'-1} \sum_{\substack{w \in \VV_k^+(j) \\ D_{ij}^{zw} \leq h_{ij}}}  \abs{u_{k,w}(x)}^2 \leq
  (d+1) \, \Clevv \sum_{k=1}^{\ell-2} \sum_{\substack{w \in \VV_k^+(j) \\ D_{ij}^{zw} \leq h_{ij}}} \abs{u_{k,w}(x)}^2.
\end{align*}
Taking the integral over $\Omega$ and the square root, multiplying by the mesh-size, as well as applying Lemma~\ref{lemma:local_estimate}, results in 
\begin{align}\label{eq:l2_generation_estimate}
  \begin{split}
  \widehat{h}_j^{-s} \, \norm[\Big]{\sum_{k=1}^{\ell'-1} \sum_{\substack{w \in \VV_k^+(j) \\ D_{ij}^{zw} \leq h_{ij}}} u_{k,w}}_{L^2(\Omega)} 
  &\leq 
  ((d+1)\Clevv)^{1/2} \Big(\sum_{k=1}^{\ell-2} \sum_{\substack{w \in \VV_k^+(j) \\ D_{ij}^{zw} \leq h_{ij}}} \widehat{h}_j^{-2s} \, \norm{u_{k,w}}_{L^2(\Omega)}^2\Big)^{1/2} \\
  &\leq 
  \Cloc \, ((d+1)\Clevv)^{1/2} \Big(\sum_{k=1}^{\ell-2} \sum_{\substack{w \in \VV_k^+(j) \\ D_{ij}^{zw} \leq h_{ij}}}
  \enorm{u_{k,w}}^2\Big)^{1/2}.
  \end{split}
\end{align}
Plugging this back into~\eqref{eq:SCS:uniformSCS} and applying a discrete Cauchy--Schwarz inequality in $(\ell',z)$, leads to 
\begin{align*}
  S_{2}(i,j) &\lesssim \delta^{j-i} \sum_{\ell'=1}^{\ell-1}\sum_{z \in \VV_{\ell'}^+(i)} \enorm{v_{\ell',z}} \, \Big(\sum_{k=1}^{\ell-2} \sum_{\substack{w \in \VV_k^+(j) \\ D_{ij}^{zw} \leq h_{ij}}}\enorm{u_{k,w}}^2\Big)^{1/2}  \\
  &\leq \delta^{j-i}
  \Big(\sum_{\ell'=1}^{\ell-1}\sum_{z \in \VV_{\ell'}^+(i)} \enorm{v_{\ell',z}}^2\Big)^{1/2} \Big(\sum_{\ell'=1}^{\ell-1}\sum_{z \in \VV_{\ell'}^+(i)}  \sum_{k=1}^{\ell-2} \sum_{\substack{w \in \VV_k^+(j) \\ D_{ij}^{zw} \leq h_{ij}}} \enorm{u_{k,w}}^2\Big)^{1/2}.
\end{align*}
We can rewrite the sum on the right-hand side as
\begin{equation*}
  \sum_{\ell'=1}^{\ell-1}\sum_{z \in \VV_{\ell'}^+(i)} \sum_{k=1}^{\ell-2} \sum_{\substack{w \in \VV_k^+(j) \\ D_{ij}^{zw} \leq h_{ij}}} \enorm{u_{k,w}}^2 
  = 
  \sum_{k=1}^{\ell-2} \sum_{w \in \VV_k^+(j)} \sum_{\ell'=1}^{\ell-1} \sum_{\substack{z \in \VV_{\ell'}^+(i) \\ D_{ji}^{wz} \leq h_{ij}}} \enorm{u_{k,w}}^2.
\end{equation*}
To estimate the inner two sums, it is sufficient to show that the number of pairs $(\ell',z)$ with $z \in \VV_{\ell'}^+(i)$ and $D_{ji}^{wz} \leq h_{ij}$ is uniformly bounded.
From \cite[Lemma 3.1]{wc2006} and \eqref{eq:SCS:level_bound}, we already know that, for a fixed vertex $z$, only a uniformly bounded number of levels $\ell'$ satisfy $z \in \VV_{\ell'}^+(i)$.
As $i \leq j$, there holds $D_{ji}^{wz} \leq h_{ji} = \widehat{h}_i.$
Hence by~\eqref{eq:SCS:distance_equivalence}, this implies $\abs{z-w} \leq C \, \widehat{h}_i$, where $C$ depends only on $\gamma$-shape regularity.
Furthermore, for all $z\neq z' \in \widehat{\VV}_i$, there holds $\abs{z-z'} \gtrsim \widehat{h}_i.$
Therefore, analogous to step 1, we can infer
\begin{equation}\label{eq:SCS:squere_packing}
  \# \set{z \in \widehat{\VV}_i \given \abs{z-w} \leq C \, \widehat{h}_i} \lesssim \Big(\frac{\widehat{h}_i}{\widehat{h}_i}\Big)^d \eqqcolon \const{C}{sq}
\end{equation}
and consequently
\begin{equation*}
  \sum_{k=1}^{\ell-2} \sum_{w \in \VV_k^+(j)} \sum_{\ell'=1}^{\ell-1} \sum_{\substack{z \in \VV_{\ell'}^+(i) \\ D_{ji}^{wz} \leq h_{ij}}} \enorm{u_{k,w}}^2
  \leq
  \const{C}{sq} \, \Clev \sum_{k=1}^{\ell-2} \sum_{w \in \VV_k^+(j)} \enorm{u_{k,w}}^2.
\end{equation*}
Altogether, for $i \leq j$, we have proven
\begin{equation}\label{eq:SCS:near_field_result}
 S_2(i,j)
  \lesssim
  \delta^{j-i} \Big(\sum_{\ell'=1}^{\ell-1}\sum_{z \in \VV_{\ell'}^+(i)} \enorm{v_{\ell',z}}^2\Big)^{1/2} \Big(\sum_{k=1}^{\ell-2}\sum_{w \in \VV_k^+(j)} \enorm{u_{k,w}}^2\Big)^{1/2}.
\end{equation}

\textbf{Step 4 (generation estimate):}
Combining~\eqref{eq:SCS:far_field_result} and~\eqref{eq:SCS:near_field_result}, it follows that
\begin{equation}\label{eq:SCS:decomp}
  \sum_{\ell'=1}^{\ell-1}\sum_{k=1}^{\ell'-1}\sum_{z \in \VV_{\ell'}^+}  \sum_{w \in \VV_k^+} \!\! a(v_{\ell',z}, u_{k,w})  
  \lesssim
  \sum_{i=0}^M\sum_{j=0}^M \widetilde{\delta}^{\abs{i-j}} \Big(\sum_{\ell'=1}^{\ell-1}\sum_{z \in \VV_{\ell'}^+(i)} \!\! \enorm{v_{\ell',z}}^2\Big)^{1/2} 
  \Big(\sum_{k=1}^{\ell-2}\sum_{w \in \VV_k^+(j)} \!\! \enorm{u_{k,w}}^2\Big)^{1/2} \!\!\!\!\!\!, 
\end{equation}
where $\widetilde{\delta} \coloneqq \max \set{\delta, 2^{-s/d}}$.
Let $(x_{i})_{i \in M}$ and $(y_{j})_{j \in M}$ be two sequences in $\R_{\geq 0}$ and $0 < q < 1$. 
Then, there holds the estimate
\begin{equation}\label{eq:SCS:ana1}
  \sum_{i=0}^M \sum_{j=0}^M q^{|i-j|} \, x_i y_j \leq \frac{2}{1-q} \Big(\sum_{i=0}^{M} x_i^2\Big)^{1/2} \Big(\sum_{j=0}^M y_j^2\Big)^{1/2}.
\end{equation}
Showing this is straightforward; nevertheless we refer to step 1 of~\cite[Proposition~5.7]{IMPS2024} for a complete proof.
Applying~\eqref{eq:SCS:ana1} to~\eqref{eq:SCS:decomp} with $q = \widetilde{\delta}$, $x_i^2 = \sum_{\ell'=1}^{\ell-1}\sum_{z \in \VV_{\ell'}^+(i)} \enorm{v_{\ell',z}}^2$, and $y_j^2 = \sum_{k=1}^{\ell-2}\sum_{w \in \VV_k^+(j)} \enorm{u_{k,w}}^2$ concludes the proof.
\end{proof}


\section{Stable decomposition}\label{sec:stable_decomp}


This section follows the proof strategy from~\cite{FFPS2017}, where a stable decomposition of lowest order discrete functions is proven in the $H^{1/2}(\Gamma)$ setting.
The strategy relies on multilevel stability of the $L^2$-projections for quasi-uniform meshes (Lemma~\ref{lemma:l2_stability}) and a pointwise bound for a Scott-Zhang quasi-interpolation operator (Lemma~\ref{lemma:value_bound_sz}).
Moreover,~\cite{KMR2025} proves a $p$-robust one level decomposition (Lemma~\ref{lemma:p_stable_decomp}) for $d=2$. Combining the two results provides an $h$- and $p$-robust multilevel decomposition in 2D (Corollary~\ref{cor:hp_stable_decomp}).

\subsection{Stability of $L^2$-projections}\label{sec:l2_projections}
Let $\widehat{\QQ}_m \colon L^2(\Omega) \to \widehat{\XX}_m^1$ denote the $L^2$-projection for $m \in \N_0$. 
There holds the following multilevel stability for $\widehat{\QQ}_m$.
The proof is analogous to~\cite[Proof of Lemma~6]{FFPS2017}.
However, for completeness we also present it here. 

\begin{lemma}\label{lemma:l2_stability}
  There exists a constant $\Cmls > 0$, such that
  \begin{equation*}
    \sum_{m=0}^\infty \widehat{h}_m^{-2s} \norm{v - \widehat{\QQ}_m v}_{L^2(\Omega)}^2 \leq \Cmls \, \enorm{v}^2 \qquad \text{for all } v \in \widetilde{H}^s(\Omega).
  \end{equation*}
  The constant $\Cmls$ depends only on $\Omega$, $d$, $s$, and $\gamma$-shape regularity.
\end{lemma}

\begin{proof}
  Let $v \in \widetilde{H}^s(\Omega)$.
  Since $\lim\limits_{j \rightarrow \infty} \norm{v - \widehat{\QQ}_j v}_{L^2(\Omega)}=0$ and $((\widehat{\QQ}_j-\widehat{\QQ}_{j-1})v,(\widehat{\QQ}_i-\widehat{\QQ}_{i-1})v)_{L^2(\Omega)}=0$ for $i \neq j$, we have
  \begin{align*}
    \norm{v-\widehat{\QQ}_m v}_{L^2(\Omega)}^2 
    &=
    \lim_{M \rightarrow \infty} \norm{(\widehat{\QQ}_M-\widehat{\QQ}_{m})v}_{L^2(\Omega)}^2 
     = 
    \lim_{M \rightarrow \infty} \norm[\Big]{\sum_{j=m+1}^M (\widehat{\QQ}_j-\widehat{\QQ}_{j-1})v}_{L^2(\Omega)}^2 \\
    &=
    \lim_{M \rightarrow \infty} \sum_{j=m+1}^M \norm{(\widehat{\QQ}_j-\widehat{\QQ}_{j-1})v}_{L^2(\Omega)}^2
    =
    \sum_{j=m+1}^\infty \norm{(\widehat{\QQ}_j-\widehat{\QQ}_{j-1})v}_{L^2(\Omega)}^2.
  \end{align*}
Moreover, since $2^{m/d}\widehat{h}_m = \widehat{h}_0$ holds for all $m \in \N_0$, the geometric series formula yields
\begin{equation*}
  \sum_{m=0}^{j} \widehat{h}_m^{-2s}
  =
  \widehat{h}_0^{-2s} \sum_{m=0}^{j} 2^{2sm/d}
  =
  \widehat{h}_0^{-2s} \, \frac{2^{2s(j+1)/d}-1}{2^{2s/d}-1}
  \leq
  \frac{\widehat{h}_0^{-2s}}{2^{2s/d}-1} 2^{2s(j+1)/d} \quad \text{for all } j \in \N_0.
\end{equation*}
With the last two observations, a change in the order of summation, and $\widehat{\QQ}_{-1} \coloneqq 0$, we obtain
\begin{align*}
  \sum_{m=0}^\infty \widehat{h}_m^{-2s} &\norm{v - \widehat{\QQ}_m v}_{L^2(\Omega)}^2 
  = 
  \sum_{m=0}^\infty \widehat{h}_m^{-2s} \sum_{j=m+1}^\infty \norm{(\widehat{\QQ}_j-\widehat{\QQ}_{j-1})v}_{L^2(\Omega)}^2 \\
  &= 
  \sum_{j=1}^{\infty} \Big(\sum_{m=0}^{j-1}\widehat{h}_m^{-{2s}}\Big) \norm{(\widehat{\QQ}_j-\widehat{\QQ}_{j-1})v}_{L^2(\Omega)}^2
  \lesssim 
  \sum_{j=0}^{\infty} 2^{2sj/d} \norm{(\widehat{\QQ}_j-\widehat{\QQ}_{j-1})v}_{L^2(\Omega)}^2.
\end{align*}
Furthermore, it is shown in~\cite[Theorem~4]{Oswald1998} and~\cite[Theorem~4.32~\&~Remark~4.33]{Xu1997} 
that, for $s \in [0,3/2)$,
\begin{equation*}
  \norm{v}_{\widetilde{H}^s(\Omega)}^2 
  \simeq 
  \sum_{j=0}^{\infty}2^{2sj/d} \norm{(\widehat{\QQ}_j-\widehat{\QQ}_{j-1})v}_{L^2(\Omega)}^2 
  \qquad 
  \text{for all } v \in \widetilde{H}^s(\Omega).
\end{equation*}
Recall $\norm{v}_{\widetilde{H}^s(\Omega)} \simeq \enorm{v}$ for $v \in \widetilde{H}^s(\Omega)$, which implies
\begin{equation*}
  \sum_{m=0}^\infty \widehat{h}_m^{-2s} \norm{v - \widehat{\QQ}_m v}_{L^2(\Omega)}^2 
  \lesssim
  \sum_{j=0}^{\infty} 2^{2sj/d} \norm{(\widehat{\QQ}_j-\widehat{\QQ}_{j-1})v}_{L^2(\Omega)}^2 
   \simeq 
   \enorm{v}^2
\end{equation*}
and the hidden constants depend only on $\Omega$, $d$, $s$, and $\gamma$-shape regularity.
This concludes the proof.
\end{proof}

\subsection{Scott-Zhang operator}
Denote by $\varphi_{\ell}^z$ the usual piecewise affine hat-function associated with the vertex $z \in \VV_{\ell}$, i.e., $\varphi_{\ell}^z|_{T} \in \P^1(T)$ for all $T \in \TT_{\ell}$ and, using the Kronecker delta notation, $\varphi_{\ell}^z(z') = \delta_{zz'}$ for all $z' \in \VV_{\ell}$. 
We consider a variant of the Scott-Zhang quasi-interpolation operator~\cite{sz1990}, which is also used in~\cite{FFPS2017} for $H^{1/2}(\Gamma)$ without boundary conditions.
Let $z \in \VV_\ell$ and $T_{\ell}^z \in \TT_{\ell}$ with $z \in T_\ell^z$.
Then, there exists a unique $L^2$-dual basis function $\psi_{\ell}^{z}$ such that
\begin{equation*}
  \int_{T_{\ell}^z} \psi_{\ell}^{z}(x)\, \varphi_{\ell}^{z'}(x) \d{x} = \delta_{zz'} \qquad \text{for all } z' \in \VV_\ell
\end{equation*}
and
\begin{equation}\label{eq:dual_norm}
  \norm{\psi_\ell^z}_{L^\infty(T_\ell^z)} \leq \Cdual \, \abs{T_\ell^z}^{-1},
\end{equation}
where $\Cdual$ is independent of the mesh $\TT_\ell$ and the vertex $z \in \VV_\ell$.
Define the operator $\JJ_{\ell} \colon L^2(\Omega) \to \XX_\ell^1$ by
\begin{equation*}
  \JJ_\ell v \coloneqq \sum_{z \in \VV_\ell \backslash \partial \Omega} \varphi_\ell^z \int_{T_\ell^z} \psi_\ell^z(x)\, v(x) \d{x} \qquad \text{for all } v \in L^2(\Omega).
\end{equation*}
The definition immediately implies $\JJ_\ell v_\ell^1 = v_\ell^1$ for all $v_\ell^1 \in \XX_\ell^1$.
For $z \in \VV_\ell \backslash \VV_\ell^+ \subseteq \VV_{\ell-1}$, there holds by definition~\eqref{def:new_vertices} that $\omega_\ell(z) = \omega_{\ell-1}(z)$ and hence $\varphi_\ell^z = \varphi_{\ell-1}^z$.
Moreover, we can choose the element $T_\ell^z$ so that $T_\ell^z = T_{\ell-1}^z$.
This implies $\psi_\ell^z = \psi_{\ell-1}^z$ and thus 
\begin{equation}\label{eq:special_property_sz}
  (\JJ_\ell-\JJ_{\ell-1})v(z) = 0 \qquad \text{for all } z \in \VV_\ell\backslash \VV_\ell^+.
\end{equation} 
Furthermore, if $v$ is affine on $T_{\ell}^z$ and vanishes on $T_{\ell}^z\cap\partial\Omega$, then
\begin{equation}\label{eq:affine_property}
  \JJ_\ell v(z) = v(z).
\end{equation}
We need another property of the operator $\JJ_\ell$, which was also shown in~\cite[Lemma~7]{FFPS2017} for the boundary element method and translates directly to the setting of this paper with adjusted powers of the mesh-size. 

\begin{lemma}[{\cite[Lemma~7]{FFPS2017}}]\label{lemma:value_bound_sz}
  For all $v \in L^2(\Omega)$, $\ell \in \N$, and $z \in \VV_\ell$, there holds
  \begin{equation}
    \abs{(\JJ_\ell-\JJ_{\ell-1})v(z)} \leq \Csz \, h_{\ell,z}^{-d/2} \, \norm{v}_{L^2(\omega_{\ell-1}^2(\{z\}))},
  \end{equation}
  where the constant $\Csz > 0$ depends only on $\gamma$-shape regularity. \qed
\end{lemma}

Finally, the proof for the stable decomposition can be performed analogously to~\cite[Proof of Proposition 4]{FFPS2017}.

\begin{proposition}[multilevel stable decomposition]\label{prop:stable_decomp}
  Let $v_{\ell}^1 \in \XX_{\ell}^1$. Then, there exist functions $v_0^1 \in \XX_0^1$ and $v_{\ell',z}^1 \in \XX_{\ell',z}^1$ for all $\ell' = 1, \dots, \ell$ and $z \in \VV_{\ell'}^+$ such that
  \begin{equation}\label{eq:decomposition}
    v_{\ell}^1 = v_0^1 + \sum_{\ell'=1}^\ell \sum_{z \in \VV_{\ell'}^+} v_{\ell',z}^1
  \end{equation}
  and
  \begin{equation}\label{eq:stability_decomposition}
    \enorm{v_0^1}^2 + \sum_{\ell'=1}^\ell \sum_{z \in \VV_{\ell'}^+} \enorm{v_{\ell',z}^1}^2 \leq \Chsd \, \enorm{v_\ell^1}^2,
  \end{equation}
  where the constant $\Chsd \geq 1$ depends only on $\Omega$, $d$, $s$, and $\gamma$-shape regularity.
\end{proposition}

\begin{proof}
This proof is split into three steps.

\textbf{Step 1 (decomposition):}
Let $v_\ell^1 \in \XX_\ell^1$ and set $\JJ_{-1} \coloneqq 0$.
Due to  the property~\eqref{eq:special_property_sz}, there exist local functions $(v^1_{\ell',z})_{z \in \VV_{\ell'}^+}$ with $v_{\ell',z}^1 \in \XX_{\ell',z}^1$ for $\ell' \in \{1,\dots,\ell\}$ and a function $v_0^1 \in \XX_0^1$ such that 
\begin{equation}\label{eq:local_decomp}
  (\JJ_{\ell'}-\JJ_{\ell'-1})v_\ell^1 = \sum_{z \in \VV_{\ell'}^+} v_{\ell',z}^1 \quad \text{and} \quad (\JJ_{0}-\JJ_{-1})v_\ell^1 = v_0^1.
\end{equation}
Together with the projection property of $\JJ_\ell$ and the telescoping sum, it follows that
\begin{equation*}
  v_\ell^1 = \sum_{\ell'=0}^{\ell} (\JJ_{\ell'}-\JJ_{\ell'-1})v_\ell^1 = v_0^1 + \sum_{\ell'=1}^\ell \sum_{z \in \VV_{\ell'}^+} v_{\ell',z}^1.
\end{equation*}
This concludes the existence of decomposition~\eqref{eq:decomposition}.

\textbf{Step 2 (stability for initial mesh):}
The triangle inequality yields
\begin{equation*}
  \enorm{\JJ_0 v_\ell^1} \leq \sum_{z \in \VV_0 \backslash \partial \Omega} \enorm{\varphi_0^z} \abs[\Big]{\int_{T_0^z} \psi_0^z(x) v_\ell^1(x) \d{x}}.
\end{equation*}
One can show via standard interpolation techniques that 
\begin{equation}\label{eq:hat_function_norm}
  \enorm{\varphi_{\ell'}^z}^2 \simeq h_{\ell',z}^{d-2s} \qquad \text{for all } \ell' \in \N_0;
\end{equation} 
see, e.g.,~\cite[Lemma 4.3]{amt1999}.
Moreover, the Cauchy--Schwarz inequality and~\eqref{eq:dual_norm} give us
\begin{equation*}
  \abs[\Big]{\int_{T_0^z} \psi_0^z(x) v_\ell^1(x) \d{x}} 
  \leq
  \norm{\psi_0^z}_{L^\infty(\Omega)} |T_0^z|^{1/2}  \norm{v_\ell^1}_{L^2(T_0^z)} 
  \eqreff{eq:dual_norm}{\lesssim}
  |T_0^z|^{-1/2} \norm{v_\ell^1}_{L^2(T_0^z)} 
  \simeq h_{0,z}^{-d/2}  \norm{v_\ell^1}_{L^2(T_0^z)}. 
\end{equation*}
Combining these observations with finite patch overlap and Lemma~\ref{lemma:local_estimate}, we derive
\begin{equation*}
  \enorm{\JJ_0 v_\ell^1} 
  \lesssim 
  \sum_{z \in \VV_0 \backslash \partial \Omega} h_{0,z}^{d/2-s} h_{0,z}^{-d/2}  \norm{v_\ell^1}_{L^2(T_0^z)}
  \lesssim 
  \diam(\Omega)^{-s} \norm{v_\ell^1}_{L^2(\Omega)}
  \eqreff{eq:local_estimate}{\lesssim}
  \enorm{v_\ell^1},
\end{equation*}
where the hidden constant depends only on the initial mesh $\TT_0$, $\gamma$-shape regularity, $\Omega$, $d$, and $s$.

\textbf{Step 3 (stability for refined meshes):}
Let $\ell' \in \{1, \dots, \ell\}$ and $z \in \VV_{\ell'}^+$.
Due to~\eqref{eq:local_decomp}, there holds $v_{\ell',z}^1 = \varphi_{\ell'}^z ((\JJ_{\ell'}-\JJ_{\ell'-1})v_\ell^1(z))$ and together with~\eqref{eq:hat_function_norm} also
\begin{equation*}
  \enorm{v_{\ell',z}^1}^2 
  = 
  \abs{(\JJ_{\ell'}-\JJ_{\ell'-1})v_\ell^1(z)}^2 \enorm{\varphi_{\ell'}^z}^2 
  \simeq 
  h_{\ell',z}^{d-2s} \, \abs{(\JJ_{\ell'}-\JJ_{\ell'-1})v_\ell^1(z)}^2.
\end{equation*}
Recall the $L^2$-projections $\widehat{\QQ}_m$ from section~\ref{sec:l2_projections} and set $\widehat{\QQ}_m \coloneqq \widehat{\QQ}_0$ for $m \in \Z$ with $m <0$.
By~\eqref{eq:level_span}, there holds $\widehat{\QQ}_{g_{\ell',z} - \Cspan} v_\ell^1 \in \widehat{\XX}_{\max\{g_{\ell',z}-\Cspan,0\}} \subseteq \widehat{\XX}^1_{r_{\ell'-1,z}}$.
With the definition~\eqref{def:smallest_level} of $r_{\ell'-1,z}$ and $r_{\ell'-1,z} \leq r_{\ell',z}$, it follows that $(\widehat{\QQ}_{g_{\ell',z} - \Cspan}v_{\ell}^1)|_T \in \P^1(T)$ for all $T \in \TT_{\ell'-1}^2(\{z\})$ and for all $T \in \TT_{\ell'}(z)$.
This immediately yields $\widehat{\QQ}_{g_{\ell',z} - \Cspan}v_{\ell}^1$ being affine on $T_{\ell'}^z \in \TT_{\ell'}(z)$ and hence $\JJ_{\ell'}(\widehat{\QQ}_{g_{\ell',z} - \Cspan}v_{\ell}^1)(z) = (\widehat{\QQ}_{g_{\ell',z} - \Cspan}v_{\ell}^1)(z)$ using~\eqref{eq:affine_property}.
Either $z \in \VV_{\ell'-1}$ and therefore $T_{\ell'-1}^z \in \TT_{\ell'-1}(z)$ or $z \in \VV_{\ell'} \backslash \VV_{\ell'-1}$ and $z$ lies on the edge spanned by two vertices $z_1,z_2 \in \VV_{\ell'-1}$ with $T_{\ell'-1}^{z_1}, T_{\ell'-1}^{z_2} \in \TT_{\ell'-1}^2(\{z\})$.
For any function $v \in L^2(\Omega)$, there holds
\begin{equation*}
  \JJ_{\ell'-1}v(z) = \varphi_{\ell'-1}^{z_1}(z) \int_{T_{\ell'-1}^{z_1}} \psi_{\ell'-1}^{z_1}(x) v(x) \d{x} + \varphi_{\ell'-1}^{z_2}(z) \int_{T_{\ell'-1}^{z_2}} \psi_{\ell'-1}^{z_2}(x) v(x) \d{x},
\end{equation*}
where $z_1 = z$ and $\varphi_{\ell'-1}^{z_2} = 0$ in the first case.
This means that the Scott-Zhang projection $\JJ_{\ell'-1}$ also preserves the value of $\widehat{\QQ}_{g_{\ell',z} - \Cspan}v_{\ell}^1$ at $z$.
Thus with Lemma~\ref{lemma:value_bound_sz}, we obtain
\begin{align*}
h_{\ell',z}^{d-2s} \, \abs{(\JJ_{\ell'}-\JJ_{\ell'-1})v_\ell^1(z)}^2 
&=
h_{\ell',z}^{d-2s} \, \abs{(\JJ_{\ell'}-\JJ_{\ell'-1})(v_\ell^1-\widehat{\QQ}_{g_{\ell',z} - \Cspan}v_{\ell}^1)(z)}^2 \\
&\lesssim 
h_{\ell',z}^{-2s} \norm{v_\ell^1-\widehat{\QQ}_{g_{\ell',z} - \Cspan}v_{\ell}^1}^2_{L^2(\omega_{\ell'-1}^2(\{z\}))}.
\end{align*}
Using $h_{\ell',z} \simeq \widehat{h}_{g_{\ell',z}}$, which holds because of $\gamma-$shape regularity, and~\eqref{eq:uniform_patch}, we get
\begin{align*}
  \sum_{\ell'=1}^{\ell} \sum_{z \in \VV_{\ell'}^+} \enorm{v_{\ell',z}^1}^2 
  &\lesssim
  \sum_{\ell'=1}^{\ell}\sum_{z \in \VV_{\ell'}^+} h_{\ell,z}^{-2s} \norm{v_\ell^1-\widehat{\QQ}_{g_{\ell',z} - \Cspan}v_{\ell}^1}^2_{L^2(\omega_{\ell'-1}^2(\{z\}))} \\
  &\simeq
  \sum_{m=0}^{\infty}\sum_{\ell'=1}^{\ell}\sum_{z \in \VV_{\ell'}^+(m)} \widehat{h}_m^{-2s} \norm{v_\ell^1-\widehat{\QQ}_{m - \Cspan}v_{\ell}^1}^2_{L^2(\omega_{\ell'-1}^2(\{z\}))} \\
  &\eqreff{eq:uniform_patch}{\leq} 
  \sum_{m=0}^{\infty}\sum_{\ell'=1}^{\ell}\sum_{z \in \VV_{\ell'}^+(m)} \widehat{h}_m^{-2s} \norm{v_\ell^1-\widehat{\QQ}_{m - \Cspan}v_{\ell}^1}^2_{L^2(\widehat{\omega}_{m}^n(z))}.
\end{align*}
Recall the set of levels $\LL_1^{\ell}(z,m)$ and the uniform bound from~\eqref{eq:SCS:level_bound}. 
Together with finite patch overlap, we further estimate
\begin{align*}
  \sum_{m=0}^{\infty}\sum_{\ell'=1}^{\ell}&\sum_{z \in \VV_{\ell'}^+(m)} \!\!\!\! \widehat{h}_m^{-2s} \norm{v_\ell^1-\widehat{\QQ}_{m - \Cspan}v_{\ell}^1}^2_{L^2(\widehat{\omega}_{m}^n(z))} 
  = 
  \sum_{m=0}^{\infty}\sum_{z \in \widehat{\VV}_m} \sum_{\ell' \in \LL_1^\ell(z,m)} \!\!\!\! \widehat{h}_m^{-2s} \norm{v_\ell^1-\widehat{\QQ}_{m - \Cspan}v_{\ell}^1}^2_{L^2(\widehat{\omega}_{m}^n(z))} \\
  &\leq 
  \Clev \sum_{m=0}^{\infty}\sum_{z \in \widehat{\VV}_m} \widehat{h}_m^{-2s} \norm{v_\ell^1-\widehat{\QQ}_{m - \Cspan}v_{\ell}^1}^2_{L^2(\widehat{\omega}_{m}^n(z))}
  \lesssim
  \sum_{m=0}^{\infty} \widehat{h}_m^{-2s} \norm{v_\ell^1-\widehat{\QQ}_{m - \Cspan}v_{\ell}^1}^2_{L^2(\Omega)}.
\end{align*}
We still need to treat the mismatch between the mesh-size function $\widehat{h}_m$ and the projection $\widehat{\QQ}_{m-\Cspan}$.
As $2^{N/d}\widehat{h}_{m+N} = \widehat{h}_m$ for any $N \in \N_0$ and $\widehat{\QQ}_m = \widehat{\QQ}_0$ for $m <0$, it follows that 
\begin{align*}
 \sum_{m=0}^{\infty} \widehat{h}_m^{-2s} \norm{v_\ell^1-&\widehat{\QQ}_{m - \Cspan}v_{\ell}^1}^2_{L^2(\Omega)}
 = 
\norm{v_\ell^1-\widehat{\QQ}_{0}v_{\ell}^1}^2_{L^2(\Omega)} \sum_{m=0}^{\Cspan-1} \widehat{h}_m^{-2s} + \sum_{m=0}^{\infty} \widehat{h}_{m+\Cspan}^{-2s}\norm{v_\ell^1-\widehat{\QQ}_{m}v_{\ell}^1}^2_{L^2(\Omega)} \\
&=
\widehat{h}_0^{-2s}\norm{v_\ell^1-\widehat{\QQ}_{0}v_{\ell}^1}^2_{L^2(\Omega)} \sum_{m=0}^{\Cspan-1}2^{2sm/d} + 2^{2s\Cspan/d} \sum_{m=0}^\infty \widehat{h}_m^{-2s} \norm{v_\ell^1-\widehat{\QQ}_{m}v_{\ell}^1}^2_{L^2(\Omega)} \\
&\lesssim  \sum_{m=0}^\infty \widehat{h}_m^{-2s} \norm{v_\ell^1-\widehat{\QQ}_{m}v_{\ell}^1}^2_{L^2(\Omega)}. 
\end{align*}
Lastly, stability of the $L^2$-projection from Lemma~\ref{lemma:l2_stability} yields 
\begin{equation*}
\sum_{\ell'=1}^{\ell} \sum_{z \in \VV_{\ell'}^+} \enorm{v_{\ell',z}^1}^2 
  \lesssim
  \sum_{m=0}^\infty \widehat{h}_m^{-2s} \norm{v_\ell^1-\widehat{\QQ}_{m}v_{\ell}^1}^2_{L^2(\Omega)}
  \lesssim \enorm{v_\ell^1}^2.
\end{equation*}
Combining steps~2 and~3 concludes the proof of stability~\eqref{eq:stability_decomposition}. 
\end{proof}

\subsection{Polynomial degree robustness}

For $d=2$, the following $p$-robust stable one level decomposition is proven in~\cite{KMR2025}.
We also refer to~\cite{KMR2020} for a stable decomposition for the boundary element method in 3D.

\begin{lemma}[one level $p$-robust decomposition]\label{lemma:p_stable_decomp}
Let $p \in \N$.
For every $v_\ell \in \XX_{\ell}^p$ there exist functions $v_{\ell}^1 \in \XX_{\ell}^1$ and $v_{\ell,z}^p \in \XX_{\ell,z}^p$ for all $z \in \VV_{\ell}$ such that
\begin{equation*}
  v_\ell = v_{\ell}^1 + \sum_{z \in \VV_\ell} v_{\ell,z}^p
\end{equation*}
and 
\begin{equation}\label{eq:one_level_stability}
  \enorm{v_{\ell}^1}^2 + \sum_{z \in \VV_\ell} \enorm{v_{\ell,z}^p}^2 \leq 
  \Cpsd
  \enorm{v_\ell}^2, 
\end{equation}
where $\Cpsd \geq 1$ depends only on $s$, $\Omega$, and $\gamma$-shape regularity.
\end{lemma}

\begin{proof}
  Let $\EE_{\Omega}$ denote the set of interior edges of the triangulation $\TT_\ell$.
  Using \cite[Lemma~4.5]{KMR2025}, $v_\ell$ admits a decomposition
  \begin{equation*}
    v_\ell 
    = 
    v_\ell^1 + \sum_{z \in \VV_\ell \backslash \partial \Omega} v_{\ell,z} + \sum_{E \in \EE_{\Omega}} v_{\ell,E} + \sum_{T \in \TT_\ell} v_{\ell,T},
  \end{equation*}
  where $v_{\ell,z} \in \XX_{\ell,z}^p$ and $v_{\ell,T}, v_{\ell,E}\in \XX_\ell^p$ with $\supp v_{\ell,E} \subseteq \overline{\omega_\ell(E)}$ and $\supp v_{\ell,T} \subseteq T$.
By~\cite[Lemma~4.5]{KMR2025}, this decomposition satisfies a series of stability estimates; however, the stability is not given directly in the energy norm to obtain~\eqref{eq:one_level_stability}.
Indeed,~\cite[Lemma~4.5]{KMR2025} provides stability of the global lowest-order, vertex, edge, and volume contributions by the full contribution in the interpolation norm $\norm{v_\ell}_{[L^2(\Omega),H_0^1(\Omega)]_s}$.
Thus, it remains to treat the left-hand side contributions that appear in different norms.
First, the global contribution is already estimated in the interpolation norm.
For the vertex and edge contributions, combining~\cite[Lemma~5.1(ii)]{KMR2025} with~\cite[Lemma~5.6]{KMR2025} yields stability in the interpolation norm.
Finally, the volume contributions can be estimated in the interpolation norm by again applying~\cite[Lemma~5.1(ii)]{KMR2025} and~\cite[Lemma~5.5]{KMR2025}.
To conclude, recall that on Lipschitz domains, the energy norm $\enorm{\cdot}$ is equivalent to the interpolation norm $\norm{\cdot}_{[L^2(\Omega),H_0^1(\Omega)]_s}$; see, e.g.,~\cite[Lemma~4.11]{CHM2015}.
Hence, it follows that
\begin{equation}\label{eq:stable_decomp_melenk}
    \enorm{v_\ell^1}^2 + \sum_{z \in \VV_\ell \backslash \partial \Omega}\enorm{v_{\ell,z}}^2 + \sum_{E \in \EE_{\Omega}} \enorm{v_{\ell,E}}^2 + \sum_{T \in \TT_\ell} \enorm{v_{\ell,T}}^2 
    \lesssim
    \enorm{v_\ell}^2,
  \end{equation}
  where the hidden constant depends only on $s$, $\Omega$, and $\gamma$-shape regularity.
  For each $v_{\ell,E}$ and $v_{\ell,T}$, choose vertices $z_E, z_T \in \VV_\ell$ such that $z_E \in E$ and $z_T \in T$, respectively. 
  Then, we can define 
  \begin{equation*}
    v_{\ell,z}^p \coloneqq v_{\ell,z} + \sum_{\substack{E \in \EE_\Omega \\ z_E = z}} v_{\ell,E} + \sum_{\substack{T \in \TT_\ell \\ z_T = z}} v_{\ell,T}.
  \end{equation*}
  Hence, the triangle inequality and $\gamma$-shape regularity together with~\eqref{eq:stable_decomp_melenk} conclude the proof.
\end{proof}

Combining the last result with Proposition~\ref{prop:stable_decomp} immediately gives the desired $h$- and $p$-robust decomposition.

\begin{corollary}[$h$- and $p$-robust stable decomposition]\label{cor:hp_stable_decomp}
Let $p \in \N $, $d=2$, and $v_{\ell} \in \XX_\ell^p$.
Then, there exist functions $v_0 \in \XX_0^1$, $v_{\ell',z} \in \XX_{\ell',z}^1$ for $\ell' \in \{1,\dots, \ell-1\}$ and $z \in \VV_{\ell'}^+$, and $v_{\ell,z} \in \XX_{\ell,z}^p$ for all $z \in \VV_\ell$ such that
\begin{equation*}
  v_\ell = v_0 + \sum_{\ell'=1}^{\ell-1} \sum_{z \in \VV_{\ell'}^+} v_{\ell',z} + \sum_{z \in \VV_{\ell}} v_{\ell,z}
\end{equation*}
and
\begin{equation*}
  \enorm{v_0}^2 + \sum_{\ell'=1}^{\ell-1} \sum_{z \in \VV_{\ell'}^+} \enorm{v_{\ell',z}}^2 + \sum_{z \in \VV_{\ell}} \enorm{v_{\ell,z}}^2 
  \leq 
  \Csd 
  \enorm{v_\ell}^2,
\end{equation*}
where $\Csd \geq 1$ depends only on $\Omega$, $s$, and $\gamma$-shape regularity. \qed
\end{corollary}

\section{Numerical examples}\label{sec:num}

All numerical examples are calculated for the following test case:
\begin{itemize}
\item L-shaped domain $\Omega := (-1,1)^2 \backslash [0,1] \times [-1,0]$, right-hand side $f = 1$.
\end{itemize}

For the matrix assembly of the stiffness matrix for the integral fractional Laplacian, we use the FEM code from \cite{ABB2017}, which supports first order piecewise polynomials $p=1$ only. The weighted residual type error estimator is taken from \cite{FMP2021} and refinement is done using 2D NVB. 
The multigrid solver is taken from the software package MooAFEM\footnote{https://gitlab.tuwien.ac.at/asc/praetorius/software/mooafem}, \cite{MooAFEM}. Calculations are done with MATLAB R2024b on an AMD Ryzen 7 5700G with 8 cores and 32GB memory.

\subsection{Convergence of AFEM with multigrid solve}
In our first numerical example, employing Algorithm~\ref{alg:afem}, we consider the convergence of the residual error estimator, where the \texttt{SOLVE} module employs the iterative multigrid. As a stopping criterion, we take $\enorm{u_\ell^k- u_\ell^{k-1}} \leq \mu \, \eta_\ell(u_\ell^0)$ with $\mu = 0.1$, i.e., the evaluation of the estimator for the initial iterate in the solver step. The evaluation of the estimator, see \cite{FMP2021}, dominates the computational cost and thus taking the estimator for the current iterate $\eta_\ell(u_\ell^k)$ as the stopping criterion would require a recomputation of the expensive estimator after each V-cycle.    

\begin{figure}[ht]
\includegraphics{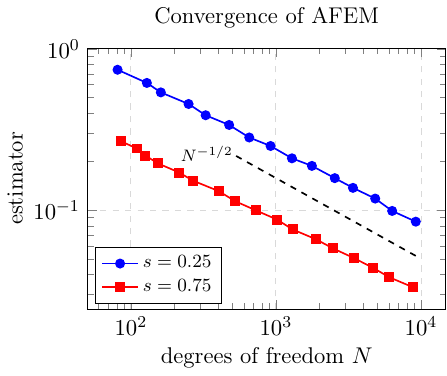}
\label{fig:convergence}
\caption{Optimal convergence of the weighted residual error estimator on adaptive NVB meshes with multigrid solver.}
\end{figure}

Figure~\ref{fig:convergence} reproduces the optimal convergence behaviour of order $N^{-1/2}$, where $N$ denotes the degrees of freedom of the FEM space. This is consistent with \cite{FMP2021}, where the same behaviour is observed for AFEM with a direct algebraic solver.  

\subsection{Comparison of solver timings}

In this experiment, for a fixed, adaptively generated mesh hierarchy on the L-shaped domain, we consider only the module \texttt{SOLVE}. In the following, we compare the direct solver using the MATLAB backslash operator with our V-cycle multigrid with optimal step-size as an iterative solver. 
We stop the multigrid once the error between consecutive iterations is below a given tolerance, chosen as $10^{-7}$ as well as using the stopping criterion $\enorm{u_\ell^k- u_\ell^{k-1}} \leq \mu \, \eta_\ell(u_\ell^0)$ with $\mu = 0.1$.

\begin{figure}[ht]
\includegraphics{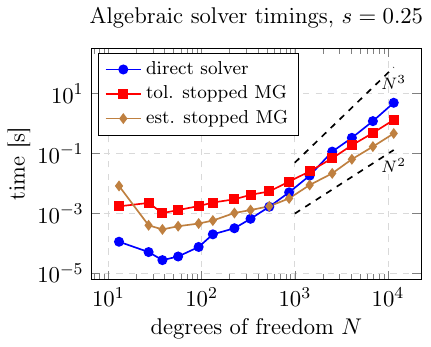}
\includegraphics{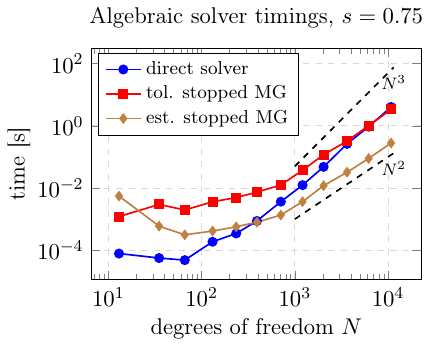}
\label{fig:timings}
\caption{Times in seconds to solve the algebraic linear system for a fixed, adaptively generated hierarchy. Left: $s=0.25,\theta = 0.3$; right: $s=0.75, \theta = 0.5$.}
\end{figure}

Figure~\ref{fig:timings} compares only the times needed to solve the algebraic linear system to calculate the FEM solution $u_L$ on the level $L$ in the fixed mesh hierarchy. In order to observe a meaningful difference between the solvers, matrix assembly is
not included in the timings. As expected, the multigrid solver behaves like $\mathcal{O}(k_{\max}N^2)$ for the fully populated stiffness matrix $A \in \R^{N\times N}$ with $k_{\max}$ being the number of V-cycles to reach the tolerance, while the direct solver scales like $\mathcal{O}(N^3)$. The estimator stopped multigrid needs significantly fewer iterations (in this example only 2 V-cycles) than the tolerance stopped multigrid.

\subsection{Contraction of the multigrid solver}

In this experiment, we compute the contraction factor $\qmg$ from \eqref{eq:contraction} for different fractional parameters $s$ on fixed, adaptively generated hierarchies with $\theta = 0.3$, where the finest mesh for each $s$ has around 5000 degrees of freedom. The multigrid is stopped once the tolerance $10^{-7}$ is reached. 

\begin{figure}[ht]
\includegraphics{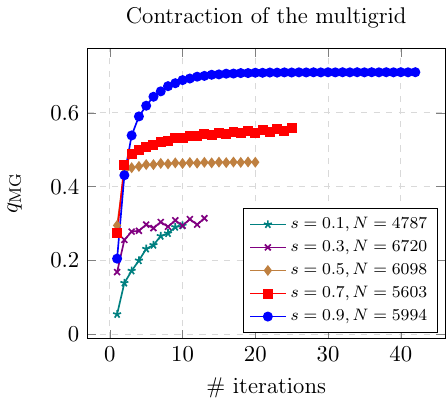}
\includegraphics{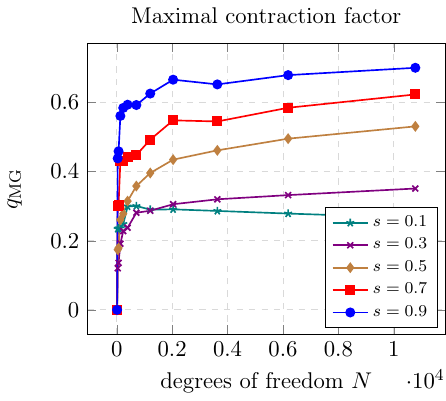}
\label{fig:contraction}
\caption{Contraction of the algebraic multigrid solver for different fractional parameters $s$. Left: contraction factor vs. number of iterations, where the finest mesh has fixed $N \sim 5000$ degrees of freedom; Right: maximal contraction factor over the mesh hierarchy vs. number of degrees of freedom.}
\end{figure}

Figure~\ref{fig:contraction} supports our theoretical result of a multigrid with a uniformly bounded contraction factor with respect to the number of iterations and degrees of freedom. Moreover, we observe a strong $s$ dependence in the contraction factor, which is not tracked in our theoretical analysis.

\subsection{Step size choice in the multigrid solver}

\begin{figure}[ht]
\includegraphics{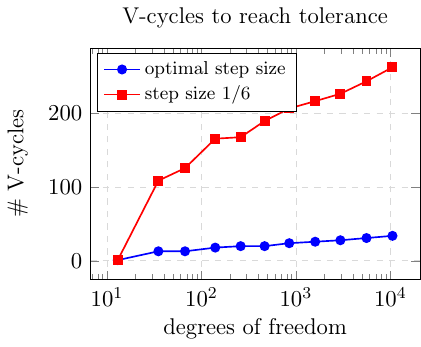}
\includegraphics{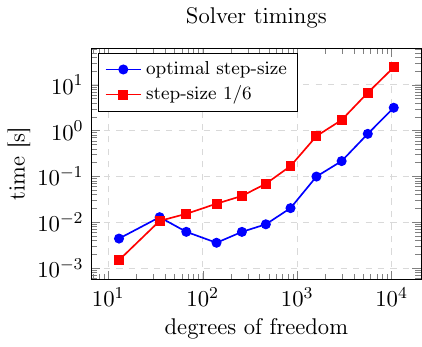}
\label{fig:step_sizes}
\caption{Step-size choice for the multigrid. Left: number of V-cycles needed to reach the tolerance of $10^{-7}$; right: solver timings for different step-size choices.}
\end{figure}

For the parameter choices $s=0.75, \theta = 0.5$, we compare the effect of different step size choices in the multigrid. In particular, in Figure~\ref{fig:step_sizes}, we compare solver times and the number of V-cycles to reach a residual tolerance of $10^{-7}$ for the optimal step size choice, which induces a bit of a computational overhead, with a fixed step size choice of $1/6 =(2(d+1))^{-1}$, for $d=2$. For our problem, clearly the optimal step size is the preferred choice.

\section*{Declaration of generative AI and AI-assisted technologies in the manuscript preparation process}

In compliance with the Leiden Declaration on Artificial Intelligence and Mathematics, the authors take full responsibility for the correctness of the mathematical arguments and results, the accuracy and completeness of the references, and the manuscript as a whole.
During the development of this work, the authors used Claude (Sonnet~4.6, 5, and~5.5) to assist
with intermediate steps, simplifications, and literature searches, and in particular with the
estimates~\eqref{eq:SCS:annuli} and~\eqref{eq:SCS:squere_packing} in the proof of
Proposition~\ref{prop:SCS}. 
For this part of the argument, Claude suggested using arguments from packing problems.

For the calculations in the numerical examples, OpenAI Codex GPT-5.6 Astra and GPT-6.1 Sol have been used to implement the crossover of the used software packages and to optimize the matrix assembly and calculation of the weighted residual error estimator both in terms of accuracy and computational times.

 After completion of the full draft of this manuscript, the authors used Claude to assist with proofreading and language editing prior to submission.

\printbibliography

\end{document}